\documentclass{article}
\usepackage{fullpage}

\usepackage[parfill]{parskip}

\usepackage[utf8]{inputenc}
\usepackage[T1]{fontenc}

\usepackage{microtype}
\microtypecontext{spacing=nonfrench}

\usepackage{amsmath,amsfonts,amssymb,mathtools,amsthm}
\usepackage{stmaryrd}
\usepackage{siunitx}
\usepackage{enumitem}
\usepackage{booktabs,multirow,multicol,arydshln}

\usepackage{xcolor}
\usepackage{hyperref}
\usepackage[noabbrev,capitalise]{cleveref}
\colorlet{inlinkcolor}{green!50!black}
\colorlet{exlinkcolor}{red!50!black}
\hypersetup{
colorlinks = true,
allcolors = inlinkcolor,
urlcolor = exlinkcolor,
}
\crefname{equation}{}{}

\usepackage{tikz}
\usepackage{caption,subcaption}
\usepackage{pgfplots}
\pgfplotsset{compat=1.16}
\usepgfplotslibrary{groupplots}

\usepackage{natbib}

\DeclareMathOperator{\logdet}{logdet}
\DeclareMathOperator{\rtdet}{rtdet}
\DeclareMathOperator{\sqr}{sqr}
\DeclareMathOperator{\spec}{spec}
\DeclareMathOperator{\lspec}{\ell_{\spec}}
\DeclareMathOperator{\linf}{\ell_{\infty}}
\DeclareMathOperator{\ltwo}{\ell_2}
\DeclareMathOperator{\geom}{geo}

\DeclareMathOperator{\wsos}{SOS}
\DeclareMathOperator{\matwsos}{matSOS}

\newcommand{\bbR}{\mathbb{R}}
\newcommand{\bbS}{\mathbb{S}}
\newcommand{\K}{\mathcal{K}}
\newcommand{\knn}{\K_{\geq}}
\newcommand{\kpsd}{\K_{\succeq}}
\newcommand{\klinf}{\K_{\linf}}
\newcommand{\kltwo}{\K_{\ltwo}}
\newcommand{\ksqr}{\K_{\sqr}}
\newcommand{\klog}{\K_{\log}}
\newcommand{\klogdet}{\K_{\logdet}}
\newcommand{\krtdet}{\K_{\rtdet}}
\newcommand{\kgeom}{\K_{\geom}}

\newcommand{\klspec}{\K_{\lspec}}
\newcommand{\kwsos}{\K_{\wsos}}
\newcommand{\kmatwsos}{\K_{\matwsos}}

\DeclareMathOperator{\cl}{cl}
\DeclareMathOperator{\diag}{diag}
\DeclareMathOperator{\Diag}{Diag}
\DeclareMathOperator{\tr}{tr}
\DeclareMathOperator{\intr}{int}
\DeclareMathOperator{\mat}{mat}
\DeclareMathOperator{\vect}{vec}
\DeclareMathOperator{\sdim}{sd}
\newcommand{\bigO}{\mathcal{O}}

\newcommand{\norm}[1]{\lVert #1 \rVert}

\newcommand{\tsum}[1]{{\textstyle \sum_{#1}}}
\newcommand{\tprod}[1]{{\textstyle \prod_{#1}}}
\newcommand{\iin}[1]{\llbracket #1 \rrbracket}

\usepackage{authblk}

\begin{document}

\title{Solving natural conic formulations with Hypatia.jl}


\author[1]{Chris Coey}
\author[1]{Lea Kapelevich}
\author[2]{Juan Pablo Vielma}
\affil[1]{Operations Research Center, MIT, Cambridge, MA}
\affil[2]{Google Research and MIT Sloan School of Management, Cambridge, MA}

\maketitle

\begin{abstract}
Many convex optimization problems can be represented through conic \emph{extended formulations} with auxiliary variables and constraints using only the small number of \emph{standard cones} recognized by advanced conic solvers such as MOSEK 9.
Such extended formulations are often significantly larger and more complex than equivalent conic \emph{natural formulations}, which can use a much broader class of \emph{exotic cones}.
We define an exotic cone as a proper cone for which we can implement tractable logarithmically homogeneous self-concordant barrier oracles for either the cone or its dual cone.
In this paper we introduce Hypatia, a highly-configurable open-source conic primal-dual interior point solver with a generic interface for exotic cones.
Hypatia is written in Julia and accessible through JuMP, and currently implements around two dozen useful predefined cones (some with multiple variants).
We define some of Hypatia's exotic cones, and for conic constraints over these cones, we analyze techniques for constructing equivalent representations using the standard cones.
For optimization problems from a variety of applications, we introduce natural formulations using these exotic cones, and we show that the natural formulations are simpler and lower-dimensional than the equivalent extended formulations.
Our computational experiments demonstrate the potential advantages, especially in terms of solve time and memory usage, of solving the natural formulations with Hypatia compared to solving the extended formulations with either Hypatia or MOSEK 9.

\end{abstract}

\setcounter{tocdepth}{3}
\tableofcontents

\section{Introduction}
\label{sec:introduction}

Any convex optimization problem may be represented as a conic problem that minimizes a linear function over the intersection of an affine subspace with a Cartesian product of primitive proper cones (i.e. irreducible, closed, convex, pointed, and full-dimensional conic sets). 
An advantage of using conic form is that, under certain conditions, a conic problem has a simple and easily checkable certificate of optimality, primal infeasibility, or dual infeasibility \citep{permenter2017solving}.
Although the scope of this paper is limited to conic problems, there are other useful notions of duality that can be leveraged by convex optimization solvers (see for example DDS solver \citep{karimi2019domain,karimi2020primal}).

\subsection{Conic interior point methods}

Most conic solvers, such as CSDP \citep{borchers1999csdp}, CVXOPT \citep{andersen2011interior}, ECOS \citep{serrano2015algorithms}, MOSEK \citep{mosek2020modeling}, and SDPA \citep{yamashita2003implementation}, implement primal-dual interior point methods (PDIPMs). 
Complexity analysis of PDIPMs, which relies on properties of \emph{logarithmically homogeneous self-concordant barrier functions} (LHSCBs; defined in \citet[Sections 2.3.1 and 2.3.3]{nesterov1994interior}), shows they require fewer iterations to converge but exhibit higher per-iteration cost compared to first order conic methods (see \citet{o2016conic} on SCS solver).
Computational evidence accords with this result and demonstrates the superior numerical robustness of PDIPMs.

Historically, PDIPM solvers were based on efficient algorithms specialized for symmetric cones, in particular, the nonnegative, (rotated) second order, and positive semidefinite (PSD) cones.
However, many useful non-symmetric conic constraints (such as $u \leq \log(w)$, representable with an exponential cone) are not exactly representable with symmetric cones.
Early non-symmetric conic PDIPMs such as \citet{nesterov1996infeasible,nesterov2012towards} had several disadvantages compared to the specialized symmetric methods, for example requiring a strictly feasible initial iterate, the solution of larger linear systems, and conjugate LHSCB oracles.
To address these issues, \citet{skajaa2015homogeneous} (henceforth referred to as \emph{SY}) introduced a PDIPM that requires just a few oracles for the primal cone only: a fixed initial point in the cone's interior, a feasibility test (to determine whether a given point is in the cone's interior), and gradient and Hessian evaluations for an LHSCB.

\subsection{Natural and extended formulations}
\label{sec:intro:natext}

Although advanced conic solvers currently recognize at most only a handful of \emph{standard cones} (nonnegative, second order, rotated second order, positive semidefinite (PSD), and $3$-dimensional exponential and power cones), these cones are sufficient for representing many problems of interest \citep{lubin2016extended,mosek2020modeling}.
Modeling tools such as disciplined convex programming (DCP) packages (see CVX \citep{grant2014cvx}, CVXPY \citep{diamond2016cvxpy}, and Convex.jl \citep{udell2014convex}) and MathOptInterface's bridges \citep{legat2020mathoptinterface} are designed to facilitate transformations of convex problems into conic problems with standard cones, to enable access to powerful specialized conic solvers.
However, for many problems of interest, a representation in terms of standard cones is not the most natural or efficient conic representation.

The process of transforming a general conic problem into a conic \emph{extended formulation} (EF) that uses only standard cones often requires introducing many artificial variables, linear equalities, and/or higher-dimensional conic constraints.
For example, in our density estimation example problem in \cref{sec:testing:densityest}, these dimensions are typically orders of magnitude larger for the EFs than for the NFs.
By increasing the size and complexity of problem data, EFs can increase the computational cost of preprocessing/initialization and linear system solving at each iteration. 
If conic solvers could recognize a much larger class of \emph{exotic} cones, they could directly solve simpler, smaller conic \emph{natural formulations} (NFs).\footnote{
We note that EFs can be beneficial for accelerating outer approximation algorithms for mixed-integer conic optimization, such as the method implemented in Pajarito solver \citep{coey2020outer}. 
However, folklore says that the EF for the second order cone likely slows down the conic solver, which is why Pajarito manages the EF in the MILP outer approximation model and only solves NFs for the conic subproblems.
}
We define an exotic cone as a proper cone for which we can implement a small set of tractable (i.e.\ fast, numerically stable, analytic) oracles for a logarithmically homogeneous self-concordant barrier for the cone or for its dual cone.

In the particular context of polynomial weighted sum-of-squares (SOS) optimization, \citet{papp2019sum} illustrate the potential advantages of using NFs with SOS cones (see \cref{sec:cones:wsos}) instead of PSD cone EFs, which are much larger.
The authors describe tractable LHSCB oracles for dual SOS cones, noting that analytic oracles are not known for primal SOS cones.
They show that their SOS NF-based approach has lower theoretical time and space complexity overall compared to a standard EF-based semidefinite programming method.
After implementing SY in their MATLAB solver Alfonso \citep{papp2017homogeneous,papp2021alfonso}, the authors observe improved solve times and scaleability from solving the NFs with Alfonso compared to solving the EFs with MOSEK.

To broaden the computational argument for NFs, in \cref{sec:cones} we define a variety of exotic cone types (some of which required the development of new LHSCBs; see \citet{coey2021self}) and describe general techniques for constructing equivalent EFs of NF constraints involving these cones.
We analyze how these EF techniques necessarily increase formulation dimensions.
We also observe that the EFs are often associated with larger values of the LHSCB \emph{parameter} $\nu$, which impacts the number of iterations $\bigO(\sqrt{\nu} \log (1 / \varepsilon))$ needed in the worst case by an idealized algorithm such as SY to obtain a solution within $\varepsilon$ tolerance \citep{nesterov1997self}.
However, most performance-oriented conic solver implementations do not directly implement idealized PDIPMs, so the practical impact of the parameter $\nu$ on average performance is unclear.
Another potential advantage of using the NFs is that converting conic certificates from the space of the EF back into the more meaningful NF space can be complicated.
Furthermore, the convergence conditions used by PDIPMs can provide numerical guarantees about conic certificates, but if EF certificates are converted to NF space, the NF certificates might lack such guarantees.

\subsection{Hypatia: a generic conic solver}

We introduce our new open source generic conic PDIPM solver, Hypatia.\footnote{
Hypatia is available at \href{https://github.com/chriscoey/Hypatia.jl}{github.com/chriscoey/Hypatia.jl} under the MIT license; see \citet{coey2021hypatiab} for documentation.
}
Hypatia is written in the Julia language \citep{bezanson2017julia} and is accessible through either a flexible, low-level native interface or the open-source modeling tools JuMP \citep{dunning2017jump} and Convex.jl \citep{udell2014convex}.
Unlike Alfonso, Hypatia uses a primal-dual form (matching CVXOPT's \emph{cone LP} form \citep{andersen2011interior}) that does not force the user to introduce slack variables, and allows linear operators to be represented with Julia's sparse, dense, or structured abstract matrix types (see \cref{sec:form}).
Hypatia already supports more than two dozen predefined exotic cone types, some of which have multiple variants (for example, real and complex flavors); these cones and associated LHSCBs are listed at \citet{coey2021hypatia,coey2021hypatiab}.

A key feature of Hypatia is the generic cone interface, which allows defining new proper cones.
The interface, like that of Alfonso, requires only the implementation of the few tractable cone oracles needed by SY.
However, unlike Alfonso, defining a new cone in Hypatia makes both the cone and its dual cone simultaneously available for use in conic formulations (see \citet{coey2021performance} for details).
For many cones of interest, tractable (i.e. fast, numerically stable, and analytic) oracles are only known for either the primal cone or the dual cone but not both.
This means Hypatia is able to handle a broader class of conic formulations than SY and Alfonso, which require oracles specifically for all cones in the primal conic formulation.
For example, in our portfolio rebalancing example NF in \cref{sec:testing:portfolio}, we have both $\ell_1$ norm cone and $\ell_\infty$ norm cone constraints; we are aware of an LHSCB with analytic oracles for the $\ell_\infty$ norm cone, but not for its dual cone - the $\ell_1$ norm cone (see \cref{sec:cones:linf}).
Unlike Alfonso, Hypatia's cone interface allows optional specification of additional cone oracles (such as dual cone feasibility tests, inverse Hessians, and third order directional derivatives; see \citet{coey2021performance}), which can improve computational efficiency and numerical performance.

Hypatia's solver interface is also highly extensible.
We provide several optional interior point search and stepping procedures, described in \citet{coey2021performance}. 
For example, while Alfonso alternates between prediction and correction steps, Hypatia's default interior point stepping procedure uses a combined directions method incorporating third-order LHSCB information, inspired by techniques of \citet{andersen2011interior,dahl2021primal,domahidi2013ecos}.
Since the per-iteration bottleneck of PDIPMs such as Hypatia's algorithm tends to be solving the large structured linear system for search directions, Hypatia allows the user to choose from several predefined methods (including options for sparse or dense factorization-based solves or linear-operator-based iterative/indirect solves) or to implement their own formulation-specific procedure to leverage additional structure.
Unlike Alfonso, Hypatia allows representing and solving conic problems in any real floating point type in Julia, hence it can solve conic problems to arbitrary precision using \emph{BigFloat} types.
We do not explore these advanced algorithmic features in this paper, and for our computational experiments, we use a fixed set of default algorithmic options.

\subsection{Examples and computational experiments}

In \cref{sec:testing} we present a series of example problems from applications such as matrix completion, experiment design, and smooth density optimization.
For these examples, we describe simple NFs using the exotic cones we define in \cref{sec:cones}.
Some of these NFs are new and may be valuable to try in real-world applications.
We randomly generate NF instances of a wide variety of sizes, construct equivalent EFs using the general EF techniques in \cref{sec:cones}, and observe that the EFs are significantly larger.

Our computational experiments demonstrate significant improvements in solve time and memory overhead from solving the NFs with Hypatia compared to solving the EFs with Hypatia or MOSEK 9.
Our experience also suggests that since EFs are often larger and more complex than NFs, they can be less convenient for the modeler, and noticeably slower and more memory-intensive to construct using JuMP or Hypatia's native interface.
For many instances, we could build the NF efficiently, but we hit time or memory limits while constructing the EF.

\section{Notation}
\label{sec:notation}

For sets, $\cl$ denotes the closure and $\intr$ denotes the interior.
$\bbR$ is the scalar reals, $\bbR_{\geq}$ is the nonnegative reals, and $\bbR_{>} = \intr ( \bbR_{\geq} )$ is the positive reals, $\bbR_{\leq}$ is the nonpositive reals, and $\bbR_{<} = \intr (\bbR_{\leq})$ is the negative reals.
The set of $d$-dimensional real vectors is $\bbR^d$, and the set of $d_1$-by-$d_2$-dimensional real matrices is $\bbR^{d_1 \times d_2}$.
$\bbS^d$ is the set of symmetric matrices of side dimension $d$, $\bbS^d_{\succeq} \subset \bbS^d$ is the positive semidefinite matrices, and $\bbS^d_{\succ} = \intr \bigl( \bbS^d_{\succeq} \bigr)$ is the positive definite matrices.
For some natural number $d$, $\iin{d}$ is the index set $\{1, 2, \ldots, d\}$.

If $a, b, c, d$ are scalars, vectors, or matrices (of appropriate dimensions), the notation $\begin{bsmallmatrix} a & b \\ c & d \end{bsmallmatrix}$ usually denotes concatenation into a matrix.
For a vector or matrix $A$, the transpose is $A'$. 
$I(d)$ is the identity matrix in $\bbR^{d \times d}$.
For dimensions implied by context, $0$ may represent vectors or matrices of $0$s, and $e$ is a vector of $1$s.
$\Diag$ represents the diagonal matrix of a given vector, and $\diag$ represents the diagonal vector of a given square matrix. 
The inner product of vectors $u, w \in \bbR^d$ is $u' w = \tsum{i \in \iin{d}} u_i w_i$.
$\log$ is the natural logarithm, $\norm{\cdot}_p$ is the $\ell_p$ norm (for $p \geq 1$) of a vector, $\det$ is the determinant of a symmetric matrix, $\tr$ is the matrix trace, and $\sigma_i(\cdot)$ is the $i$th largest singular value of a matrix.

The operator $\vect$ maps $\bbR^{d_1 \times d_2}$ (matrices) to $\bbR^{d_1 d_2}$ (vectors) by stacking columns.
The inverse operator $\mat_{d_1, d_2}$ maps $\bbR^{d_1 d_2}$ to $\bbR^{d_1 \times d_2}$.
For symmetric matrices, $\vect$ maps $\bbS^d$ to $\bbR^{\sdim(d)}$, where we define $\sdim(d) \coloneqq d (d + 1) / 2$, by rescaling off-diagonal elements by $\sqrt{2}$ and stacking columns of the upper triangle. 
For example, for $S \in \bbS^3$ we have $\sdim(3) = 6$ and $\vect(S) = \bigl( S_{1,1}, \sqrt{2} S_{1,2}, S_{2,2}, \sqrt{2} S_{1,3}, \sqrt{2} S_{2,3}, S_{3,3} \bigr) \in \bbR^{\sdim(3)}$.
The inverse mapping $\mat$ from $\bbR^{\sdim(d)}$ to $\bbS^d$ is well-defined.
The linear operators $\vect$ and $\mat$ preserve inner products, so $\vect(S)' \vect(Z) = \tr(S' Z)$ for $S, Z \in \bbR^{d_1 \times d_2}$ or $S, Z \in \bbS^d$.

\section{Conic duality and standard form}
\label{sec:form}

Let $\K$ be a proper cone in $\bbR^q$, i.e. a conic subset of $\bbR^q$ that is closed, convex, pointed, and full-dimensional (see \citet{skajaa2015homogeneous}).
We call $\K$ a primitive (or irreducible) cone if it cannot be written as a Cartesian product of two or more lower-dimensional cones.
$\K^\ast \subset \bbR^q$ is the dual cone of $\K$:
\begin{equation}
\K^\ast \coloneqq \{ z \in \bbR^q : s' z \geq 0, \forall s \in \K \}.
\label{eq:Kdual}
\end{equation}
$\K^\ast$ is a primitive proper cone if and only if $\K$ is a primitive proper cone.

In \cref{sec:cones}, we introduce a subset of Hypatia's predefined primitive proper cones and their dual cones.
We use these cones to formulate NFs and EFs for our applied example problems in \cref{sec:testing}.
In this paper, we omit the LHSCBs for the cones in \cref{sec:cones}, but these can be found in \citet{coey2021hypatia,kapelevich2021sum,coey2021self}.
We note in \cref{sec:cones} that for many of these cones, computing conjugate barrier oracles requires optimization, which is slow and numerically fraught in our experience.
Fortunately, like the algorithm by \citet{skajaa2015homogeneous} (and its implementation in Alfonso), Hypatia does not use conjugate barrier oracles.

Hypatia's generic cone interface allows defining any proper cone $\K$ by specifying a small list of oracles: an initial interior point $t \in \intr(\K)$, a feasibility test for $\intr(\K)$, and gradients and Hessians of an LHSCB $f$ for $\K$.
Recall that \citet[Sections 2.3.1 and 2.3.3]{nesterov1994interior} defines an LHSCB for a proper cone.
The cone interface also allows optional specification of other oracles that can improve performance.
Once defined, the cone and its dual cone may be used in any combination with other cones recognized by Hypatia to construct the Cartesian product cone $\K$ in the primal conic form \cref{eq:prim} below.

Hypatia's primal conic form over variable $x \in \bbR^n$ is:
\begin{subequations}
\begin{align}
\textstyle\inf_x \quad c' x &:
\label{eq:prim:obj}
\\
b - A x &= 0,
\label{eq:prim:eq}
\\
h - G x &\in \K,
\label{eq:prim:K}
\end{align}
\label{eq:prim}
\end{subequations}
where $c \in \bbR^n$, $b \in \bbR^p$, and $h \in \bbR^q$ are vectors, $A : \bbR^n \to \bbR^p$ and $G : \bbR^n \to \bbR^q$ are linear maps, and $\K \subset \bbR^q$ is a Cartesian product $\K = \K_1 \times \cdots \times \K_K$ of primitive proper cones.
Henceforth we use $n, p, q$ to denote the variable, equality, and conic constraint dimensions of a conic problem in the form \cref{eq:prim}.
This primal form matches CVXOPT's form, though CVXOPT only recognizes symmetric cones \citep{vandenberghe2010cvxopt}.
Unlike the conic form used by \citet{skajaa2015homogeneous} (and Alfonso and MOSEK 9), which recognizes conic constraints of the form $x \in \K$, our form does not require the user to introduce slack variables to represent a more general constraint $h - G x \in \K$.

The corresponding conic dual problem over variable $y \in \bbR^p$ associated with \cref{eq:prim:eq}, and $z \in \bbR^q$ associated with \cref{eq:prim:K}, is:
\begin{subequations}
\begin{align}
\textstyle\sup_{y, z} \quad -b' y - h' z &:
\label{eq:dual:obj}
\\
c + A' y + G' z &= 0,
\label{eq:dual:eq}
\\
z &\in \K^\ast,
\label{eq:dual:K}
\end{align}
\label{eq:dual}
\end{subequations}
where \cref{eq:dual:eq} is associated with primal variable $x \in \bbR^n$.
Note $\K^\ast = \K_1^\ast \times \cdots \times \K_K^\ast$.

If neither the primal nor the dual is \emph{ill-posed},\footnote{
Intuitively, according to \citet[Section 7.2]{mosek2020modeling}, a conic problem is ill-posed if a small perturbation of the problem data can change the feasibility status of the problem or cause arbitrarily large perturbations to the optimal solution.
See \citet{permenter2017solving} for more details.
} 
there exists a simple conic certificate providing an easily verifiable proof of infeasibility of the primal \cref{eq:prim} or dual \cref{eq:dual} or optimality of a given primal-dual solution.
A primal improving ray $x$ is a feasible direction for the primal along which the primal objective improves (i.e. $c' x < 0$, $-A x = 0$, $-G x \in \K$), and hence it certifies dual infeasibility via the conic generalization of Farkas' lemma.
Similarly, a dual improving ray $(y, z)$ certifies primal infeasibility (i.e. $-b' y - h' z > 0$, $A' y + G' z = 0$, $z \in \K^\ast$).
Finally, a complementary solution $(x, y, z)$ satisfies the primal-dual feasibility conditions \crefrange{eq:prim:eq}{eq:prim:K} and \crefrange{eq:dual:eq}{eq:dual:K} and has equal and attained objective values $c' x = -b' y - h' z$, and hence certifies optimality of $(x, y, z)$ via conic weak duality.

\section{Cones and extended formulations}
\label{sec:cones}

Recall that we define the standard cones as those recognized by MOSEK 9, listed below.
\begin{description}
\item[Nonnegative cone.]
The self-dual nonnegative cone is $\knn = \knn^\ast \coloneqq \bbR_{\geq}$.
\item[Euclidean norm cone.]
The self-dual Euclidean norm cone (or second order cone) is the epigraph of the $\ell_2$ norm: 
\begin{equation}
\kltwo = \kltwo^\ast \coloneqq \bigl\{ 
(u, w) \in \bbR_{\geq} \times \bbR^d : u \geq \norm{w}
\bigr\}.
\end{equation}
\item[Euclidean norm-squared cone.]
The self-dual Euclidean norm-squared cone (or rotated second order cone) is the epigraph of the perspective function of $g$ for $g(w) = \frac{1}{2} \norm{w}^2$:
\begin{equation}
\ksqr = \ksqr^\ast \coloneqq \bigl\{ 
(u, v, w) \in \bbR_{\geq}^2 \times \bbR^d : 2 u v \geq \norm{w}^2 
\bigr\}.
\end{equation}
\item[Positive semidefinite cone.]
The self-dual (vectorized) positive semidefinite (PSD) cone is:
\begin{equation}
\kpsd = \kpsd^\ast \coloneqq \bigl\{ 
w \in \bbR^{\sdim(d)} : \mat(w) \in \bbS^d_{\succeq} 
\bigr\}.
\end{equation}
\item[$3$-dimensional exponential cone.]
The exponential cone in $\bbR^3$ is a special case of our logarithm cone $\klog$ defined in \cref{sec:cones:log} (let $d = 1$ in \cref{eq:log}), so any $3$-dimensional $\klog$ constraint is an exponential cone constraint.
\item[$3$-dimensional power cone.]
The power cone in $\bbR^3$ (defined in \citet[Section 4.1]{mosek2020modeling}) is a special case of Hypatia's generalized power cone (see \citet{coey2021hypatia}).
However, none of our example NFs or EFs in \cref{sec:testing} need power cones, so we omit these definitions here.
\end{description}

In \crefrange{sec:cones:linf}{sec:cones:wsos}, we define a subset of Hypatia's predefined exotic cones.
For simplicity, we refer to a particular exotic cone constraint as an NF, and an equivalent reformulation of such a constraint in terms of only standard cones as an EF.
We describe general techniques for constructing EFs for the types of exotic conic constraints we use in our example NFs in \cref{sec:testing}.
In general, an NF constraint has the form $h - G x \in \K$, but in this section we write $s \in \K$ for simplicity, since $s = h - G x$ can be substituted into the EF description.
An EF may use auxiliary variables, linear equalities, and/or conic constraints, which affect the dimensions $n$, $p$, and $q$ (respectively) of the primal conic form \cref{eq:prim}.
In \cref{tab:natvsext}, we compare the dimensions and LHSCB parameters ($\nu$) associated with equivalent NF and EF constraints; as we mention in \cref{sec:intro:natext}, these properties may affect the performance of PDIPMs.

The EFs we describe below follow best practices from DCP modeling tools such as Convex.jl \citep{udell2014convex} and descriptions such as \citet[chapter 4]{ben2001lectures}.
We use JuMP \citep{dunning2017jump} to build the NFs and EFs in \cref{sec:testing}, so we contributed several exotic cones and the EFs described in \cref{sec:cones:linf,sec:cones:geom,sec:cones:spec} to MathOptInterface's bridges \citep{legat2020mathoptinterface} to permit automated EF construction.
MathOptInterface does not currently recognize $\klog$ (for $d > 1$), $\kwsos$, and $\kmatwsos$ (or their dual cones), so we construct the EFs in \cref{sec:cones:log,sec:cones:wsos} manually using JuMP.
For some EFs with auxiliary variables and equalities, it is possible to perform eliminations to reduce dimensions slightly, but this can impact the sparsity of problem data (note that in our experiments in \cref{sec:testing}, both Hypatia and MOSEK perform preprocessing).

\begin{table}[!htb]
\centering
\caption{
Properties of NFs and EFs for the exotic conic constraints in \crefrange{sec:cones:linf}{sec:cones:wsos}.
$q$ and $\nu$ are the dimension and LHSCB parameter for the NF cone, and $\bar{q}$ and $\bar{\nu}$ are the corresponding values for the EF Cartesian product cone. 
$\bar{n}$ and $\bar{p}$ are the EF auxiliary variable and equality dimensions. 
Note $\sdim(k)$ is $\bigO(k^2)$.
}
\label{tab:natvsext}
\begin{tabular}{lcclcccc}
\toprule
NF & $q$ & $\nu$ & EF & $\bar{q}$ & $\bar{\nu}$ & $\bar{n}$ & $\bar{p}$ \\ 
\midrule
$\klinf$ & $1 + d$ & $1 + d$ & $\knn$ & $2d$ & $2d$ & 0 & 0 \\
$\klinf^\ast$ & $1 + d$ & $1 + d$ & $\knn$ & $1 + 2d$ & $1 + 2d$ & $2d$ & $d$ \\
$\klspec$ & $1 + r s$ & $1 + r$ & $\kpsd$ & $\sdim(r + s)$ & $r + s$ & 0 & 0 \\
$\klspec^\ast$ & $1 + r s$ & $1 + r$ & $\knn, \kpsd$ & $1 + \sdim(r + s)$ & $1 + r + s$ & $\sdim(r) + \sdim(s)$ & 0 \\
$\kgeom$ & $1 + d$ & $1 + d$ & $\knn, \klog$ & $2 + 3d$ & $2 + 3d$ & $1 + d$ & 0 \\
$\krtdet$ & $1 + \sdim(d)$ & $1 + d$ & $\knn, \klog, \kpsd$ & $2 + 3d + \sdim(2d)$ & $2 + 5d$ & $1 + d + \sdim(d)$ & 0 \\
$\klog$ & $2 + d$ & $2 + d$ & $\knn, \klog$ & $1 + 3d$ & $1 + 3d$ & $d$ & 0 \\
$\klogdet$ & $2 + \sdim(d)$ & $2 + d$ & $\knn, \klog, \kpsd$ & $1 + 3d + \sdim(2d)$ & $1 + 5d$ & $1 + d + \sdim(d)$ & 0 \\
$\kwsos$ & $d$ & $\tsum{l} t_l$ & $\kpsd$ & $\tsum{l} \sdim(t_l)$ & $\nu$ & $\bar{q}$ & $d$ \\
$\kwsos^\ast$ & $d$ & $\tsum{l} s_l$ & $\kpsd$ & $\tsum{l} \sdim(s_l)$ & $\nu$ & 0 & 0 \\
$\kmatwsos$ & $\sdim(t) d$ & $t \tsum{l} s_l$ & $\kpsd$ & $\tsum{l} \sdim(t s_l)$ & $\nu$ & $\bar{q}$ & $q$ \\
\bottomrule
\end{tabular}
\end{table}

\subsection{Infinity norm cone}
\label{sec:cones:linf}

The $\ell_\infty$ norm cone is the epigraph of $\ell_\infty$, and its dual cone is the $\ell_1$ norm cone:
\begin{subequations}
\begin{align}
\klinf & \coloneqq \bigl\{ 
(u, w) \in \bbR_{\geq} \times \bbR^d : 
u \geq \norm{w}_\infty 
\bigr\},
\\
\klinf^\ast & \coloneqq \bigl\{ 
(u, w) \in \bbR_{\geq} \times \bbR^d : 
u \geq \norm{w}_1 
\bigr\}.
\end{align}
\end{subequations}
For $\klinf$, we use the LHSCB from \citet[Section 7.5]{guler1996barrier}.
We are not aware of an LHSCB for $\klinf^\ast$ with similarly efficient oracles.

Our examples in \cref{sec:testing:portfolio,sec:testing:expdesign} use the following NF (left) and EF (right):
\begin{equation}
(u, w) \in \klinf \subset \bbR^{1 + d} 
\quad \Leftrightarrow \quad
(u e - w, u e + w) \in (\knn)^{2d},
\label{eq:cones:linf}
\end{equation}
and similarly, \cref{sec:testing:portfolio} uses:
\begin{equation}
(u, w) \in \klinf^\ast \subset \bbR^{1 + d} 
\quad \Leftrightarrow \quad
\exists \theta \in (\knn)^d, \exists \lambda \in (\knn)^d,
w = \theta - \lambda, u - e' (\theta + \lambda) \in \knn.
\label{eq:cones:l1}
\end{equation}

\subsection{Spectral norm cone}
\label{sec:cones:spec}

The spectral norm cone is the epigraph of the matrix spectral norm, and its dual cone is the epigraph of the matrix nuclear norm:
\begin{subequations}
\begin{align}
\K_{\lspec(r, s)} & \coloneqq \bigl\{ 
(u, w) \in \bbR_{\geq} \times \bbR^{r s} : 
u \geq \sigma_1(W) 
\bigr\},
\\
\K_{\lspec(r, s)}^\ast & \coloneqq \bigl\{ 
(u, w) \in \bbR_{\geq} \times \bbR^{r s} : 
u \geq \tsum{i \in \iin{r}} \sigma_i(W) 
\bigr\},
\end{align} 
\end{subequations}
where $W \coloneqq \mat_{r, s}(w) \in \bbR^{r \times s}$ and $r \leq s$ (this is nonrestrictive since the singular values are the same for $W$ and $W'$).
For $\klspec$ we use the LHSCB from \citet{nesterov1994interior}.
We are not aware of an LHSCB for $\klspec^\ast$ with similarly efficient oracles. 

\Cref{sec:testing:matrixcompletion} uses the EF from \citet[Section 4.2]{ben2001lectures}:
\begin{equation}
(u, w) \in \K_{\lspec(r, s)} \subset \bbR^{1 + r s} 
\quad \Leftrightarrow \quad
\begin{bmatrix}
u I(r) & W \\
W' & u I(s)
\end{bmatrix}
\in \bbS^{r + s}_{\succeq}.
\label{eq:cones:spec}
\end{equation}
\Cref{sec:testing:matrixregression} uses the EF from \citet{recht2010guaranteed}:
\begin{align}
(u, w) \in \K_{\lspec(r, s)}^\ast \subset \bbR^{1 + r s} 
\quad \Leftrightarrow \quad
\begin{split}
& \exists \theta \in \bbR^{\sdim(r)}, \exists \lambda \in \bbR^{\sdim(s)},
\begin{bmatrix}
\Theta & W \\
W' & \Lambda
\end{bmatrix}
\in \bbS^{r + s}_{\succeq},
\\
& u - ( \tr(\Theta) + \tr(\Lambda) ) / 2 \in \knn,
\end{split}
\label{eq:cones:nuc}
\end{align}
where $\Theta \coloneqq \mat(\theta) \in \bbS^r, \Lambda \coloneqq \mat(\lambda) \in \bbS^s$.

\subsection{Geometric mean cone}
\label{sec:cones:geom}

The geometric mean cone is the hypograph of the geometric mean function:
\begin{subequations}
\begin{align}
\kgeom & \coloneqq \bigl\{ 
(u, w) \in \bbR \times \bbR^d_{\geq} : 
u \leq \tprod{i \in \iin{d}} w_i^{1/d} 
\bigr\},
\label{eq:geom}
\\
\kgeom^\ast & \coloneqq \bigl\{ 
(u, w) \in \bbR_{\leq} \times \bbR^d_{\geq} : 
u \geq -d \tprod{i \in \iin{d}} w_i^{1/d} 
\bigr\}.
\end{align}
\end{subequations}
For $\kgeom$ we use the LHSCB from \citet{nesterov2006constructing}.

The example in \cref{sec:testing:matrixcompletion} uses an EF for $\kgeom$, and the root-determinant variant of the example in \cref{sec:testing:expdesign} uses an EF for $\kgeom$ indirectly through a $\krtdet$ EF (see \cref{sec:cones:rootdet}). 
We are aware of three EFs for $\kgeom$: a rotated second order cone EF (\emph{EF-sec}) from \citet[Section 3.3.1]{ben2001lectures}, a power cone EF (\emph{EF-pow}) from \citet{mosek2020modeling}, and an exponential cone EF (\emph{EF-exp}).
We contributed EF-exp to MathOptInterface as a combination of two bridges (geometric mean cone to relative entropy cone to exponential cones):
\begin{align}
(u, w) \in \kgeom \subset \bbR^{1 + d} 
\quad \Leftrightarrow \quad
\begin{split}
& \exists \theta \in \knn, \exists \lambda \in \bbR^d, 
e' \lambda \in \knn,
\\
& (\lambda_i, u + \theta, w_i) \in \klog, \forall i \in \iin{d}.
\end{split}
\label{eq:cones:geomlog}
\end{align}
EF-pow is not currently available through MathOptInterface bridges, and it has a very similar size and structure to EF-exp, so we do not describe or test it.
EF-sec uses multiple levels of variables and $3$-dimensional $\ksqr$ constraints and is complex to describe, so we refer the reader to \citet[Section 3.3.1]{ben2001lectures}.
In our empirical comparisons in \cref{sec:testing:matrixcompletion,sec:testing:expdesign}, EF-sec typically has larger variable and conic constraint dimensions but smaller barrier parameter than EF-exp.

\subsection{Root-determinant cone}
\label{sec:cones:rootdet}

The root-determinant cone is the hypograph of the root-determinant function:
\begin{subequations}
\begin{align}
\krtdet & \coloneqq \bigl\{ 
(u, w) \in \bbR^{1 + \sdim(d)} : 
W \in \bbS^d_{\succeq}, u \leq (\det(W))^{1/d} 
\bigr\},
\label{eq:rootdet}
\\
\krtdet^\ast & \coloneqq \bigl\{ 
(u, w) \in \bbR_{\leq} \times \bbR^{\sdim(d)} : 
W \in \bbS^d_{\succeq}, u \geq -d (\det(W))^{1/d} 
\bigr\},
\end{align}
\end{subequations}
where $W \coloneqq \mat(w)$.
For $\krtdet$ we propose an LHSCB with efficient oracles in \citet{coey2021self}.

\Cref{sec:testing:expdesign} uses the EF from \citet[Section 4.2]{ben2001lectures}:
\begin{align}
(u, w) \in \krtdet \subset \bbR^{1 + \sdim(d)} 
\quad \Leftrightarrow \quad
\begin{split}
\exists \theta \in \bbR^{\sdim(d)}, 
(u, \diag(\Theta)) \in \kgeom,
\\
\begin{bmatrix} 
W & \Theta \\ 
\Theta' & \Diag(\diag(\Theta)) 
\end{bmatrix} 
\in \bbS^{2d}_{\succeq},
\end{split}
\label{eq:cones:rootdet}
\end{align}
where $\Theta \coloneqq \mat(\theta) \in \bbS^d$, and the $\kgeom$ constraint is itself replaced with one of the geometric mean cone EFs described in \cref{sec:cones:geom}.

\subsection{Logarithm cone}
\label{sec:cones:log}

The logarithm cone is the hypograph of the perspective function of a sum of natural log functions:
\begin{subequations}
\begin{align}
\klog & \coloneqq \cl \bigl\{ 
(u, v, w) \in \bbR \times \bbR^{1 + d}_{>} : 
u \leq \tsum{i \in \iin{d}} v \log ( w_i / v ) 
\bigr\},
\label{eq:log}
\\
\klog^\ast & \coloneqq \cl \bigl\{ 
(u, v, w) \in \bbR_{<} \times \bbR \times \bbR^d_{>} : 
v \geq \tsum{i \in \iin{d}} u ( \log ( -w_i / u ) + 1 ) 
\bigr\}.
\end{align}
\end{subequations}
For $\klog$ we propose an LHSCB with efficient oracles in \citet{coey2021self}.

\Cref{sec:testing:densityest} uses the EF (when $d > 1$):
\begin{equation}
(u, v, w) \in \klog \subset \bbR^{2 + d} 
\quad \Leftrightarrow \quad
\exists \theta \in \bbR^d, e' \theta - u \in \knn, 
(\theta_i, 1, w_i) \in \klog, \forall i \in \iin{d}.
\label{eq:cones:log}
\end{equation}

\subsection{Log-determinant cone}
\label{sec:cones:logdet}

The log-determinant cone is the hypograph of the perspective function of the log-determinant function: 
\begin{subequations}
\begin{align}
\klogdet & \coloneqq \cl \bigl\{ 
(u, v, w) \in \bbR \times \bbR_{>} \times \bbR^{\sdim(d)} : 
W \in \bbS^d_{\succ}, 
u \leq v \logdet( W / v )
\bigr\},
\label{eq:logdet}
\\
\klogdet^\ast & \coloneqq \cl \bigl\{ 
(u, v, w) \in \bbR_{<} \times \bbR^{1 + \sdim(d)} : 
W \in \bbS^d_{\succ}, 
v \geq u ( \logdet( -W / u ) + d ) 
\bigr\},
\end{align}
\end{subequations}
where $W \coloneqq \mat(w)$.
For $\klogdet$ we propose an LHSCB with efficient oracles in \citet{coey2021self}.

\Cref{sec:testing:expdesign} adapts the root-determinant cone EF \cref{eq:cones:rootdet}:
\begin{align}
(u, v, w) \in \klogdet \subset \bbR^{2 + \sdim(d)} 
\quad \Leftrightarrow \quad
\begin{split}
\exists \theta \in \bbR^{\sdim(d)}, 
(u, v, \diag(\Theta)) \in \klog,
\\
\begin{bmatrix} 
W & \Theta \\ 
\Theta' & \Diag(\diag(\Theta)) 
\end{bmatrix} 
\in \bbS^{2d}_{\succeq},
\end{split}
\label{eq:cones:logdet}
\end{align}
where $\Theta \coloneqq \mat(\theta) \in \bbS^d$, and the $\klog$ constraint is itself replaced with the logarithm cone EF described in \cref{sec:cones:log}.

\subsection{Polynomial weighted SOS scalar and matrix cones}
\label{sec:cones:wsos}

Given a collection of matrices $P_l \in \bbR^{d \times s_l}, \forall l \in \iin{r}$ derived from basis polynomials evaluated at $d$ interpolation points as in \citet{papp2019sum}, the interpolant basis polynomial weighted SOS cone is:
\begin{subequations}
\begin{align}
\K_{\wsos(P)} & \coloneqq \bigl\{ 
w \in \bbR^d : 
\exists \Theta_l \in \bbS^{s_l}_{\succeq}, \forall l \in \iin{r}, 
w = \tsum{l \in \iin{r}} \diag \bigl( P_l \Theta_l P_l' \bigr) 
\bigr\},
\label{eq:wsos:prim}
\\
\K_{\wsos(P)}^\ast & \coloneqq \bigl\{ 
w \in \bbR^d : 
P_l' \Diag(w) P_l \in \bbS^{s_l}_{\succeq}, \forall l \in \iin{r} 
\bigr\}.
\label{eq:wsos:dual}
\end{align}
\label{eq:wsos}
\end{subequations}
These cones are useful for polynomial and moment modeling; for example, a point in $\K_{\wsos(P)}$ corresponds to a polynomial that is pointwise nonnegative on a semialgebraic domain defined by $P$.

Given a side dimension $t$ of a symmetric matrix of polynomials (for simplicity, all using the same interpolant basis), and $P_l \in \bbR^{d \times s_l}, \forall l \in \iin{r}$ defined as for $\K_{\wsos(P)}$ in \cref{sec:cones:wsos}, the interpolant basis polynomial weighted SOS matrix cone is:
\begin{subequations}
\begin{align}
\K_{\matwsos(P)} & \coloneqq \Biggl\{
\begin{aligned} 
& w \in \bbR^{\sdim(t) d} : 
\exists \Theta_l \in \bbS^{s_l t}_{\succeq}, \forall l \in \iin{r},
\\
& W_{i,j,:} = \tsum{l \in \iin{r}} \diag \bigl( P_l (\Theta_l)_{i,j} P_l' \bigr), 
\forall i, j \in \iin{t} : i \geq j
\end{aligned}
\Biggr\},
\label{eq:matwsos:prim}
\\
\K_{\matwsos(P)}^\ast & \coloneqq \bigl\{ 
w \in \bbR^{\sdim(t) d} : 
\bigl[ P_l' \Diag(W_{i,j,:}) P_l \bigr]_{i,j \in \iin{t}} 
\in \bbS^{s_l t}_{\succeq}, \forall l \in \iin{r} 
\bigr\},
\label{eq:matwsos:dual}
\end{align}
\end{subequations}
where $W_{i,j,:} \in \bbR^d$ is the contiguous slice of $w$ (scaled to account for symmetry) corresponding to the interpolant basis values in the $(i,j)$th position of the symmetric matrix,
$(S)_{i,j}$ is the $(i,j)$th block in a symmetric matrix $S$ with square blocks of equal dimensions, and $[g(W_{i,j,:})]_{i,j \in \iin{t}}$ is the symmetric matrix with square matrix $g(W_{i,j,:})$ in the $(i,j)$th block.
A point in $\K_{\matwsos(P)}$ corresponds to a polynomial matrix that is pointwise PSD on a semialgebraic domain defined by $P$.
See \citet{kapelevich2021sum} for more details.

\Citet{papp2019sum} describe an LHSCB with efficient oracles for $\K_{\wsos(P)}^\ast$, but they state that one is not known for $\K_{\wsos(P)}$.
We are not aware of a useful LHSCB for $\K_{\matwsos(P)}$; indeed, for $t = 1$, $\K_{\matwsos(P)}$ reduces to $\K_{\wsos(P)}$.
Noting that \cref{eq:matwsos:dual} implicitly constrains a linear function of $w$ to a Cartesian product of PSD cones, we can use \citet[Propositions 5.1.1 and 5.1.3]{nesterov1994interior} (with the $-\logdet$ LHSCB for $\kpsd$) to derive an LHSCB with efficient oracles for $\K_{\matwsos(P)}^\ast$.
This LHSCB reduces to the $\K_{\wsos(P)}$ LHSCB for $t = 1$.

Our examples in \cref{sec:testing:polymin,sec:testing:densityest,sec:testing:shapeconregr} use the EFs implicit in the definitions of $\kwsos^\ast$, $\kwsos$, and $\kmatwsos$ in \cref{eq:wsos:prim,eq:wsos:dual,eq:matwsos:prim}.
These EFs each use $r$ $\kpsd$ cones.

\section{Numerical examples}
\label{sec:testing}

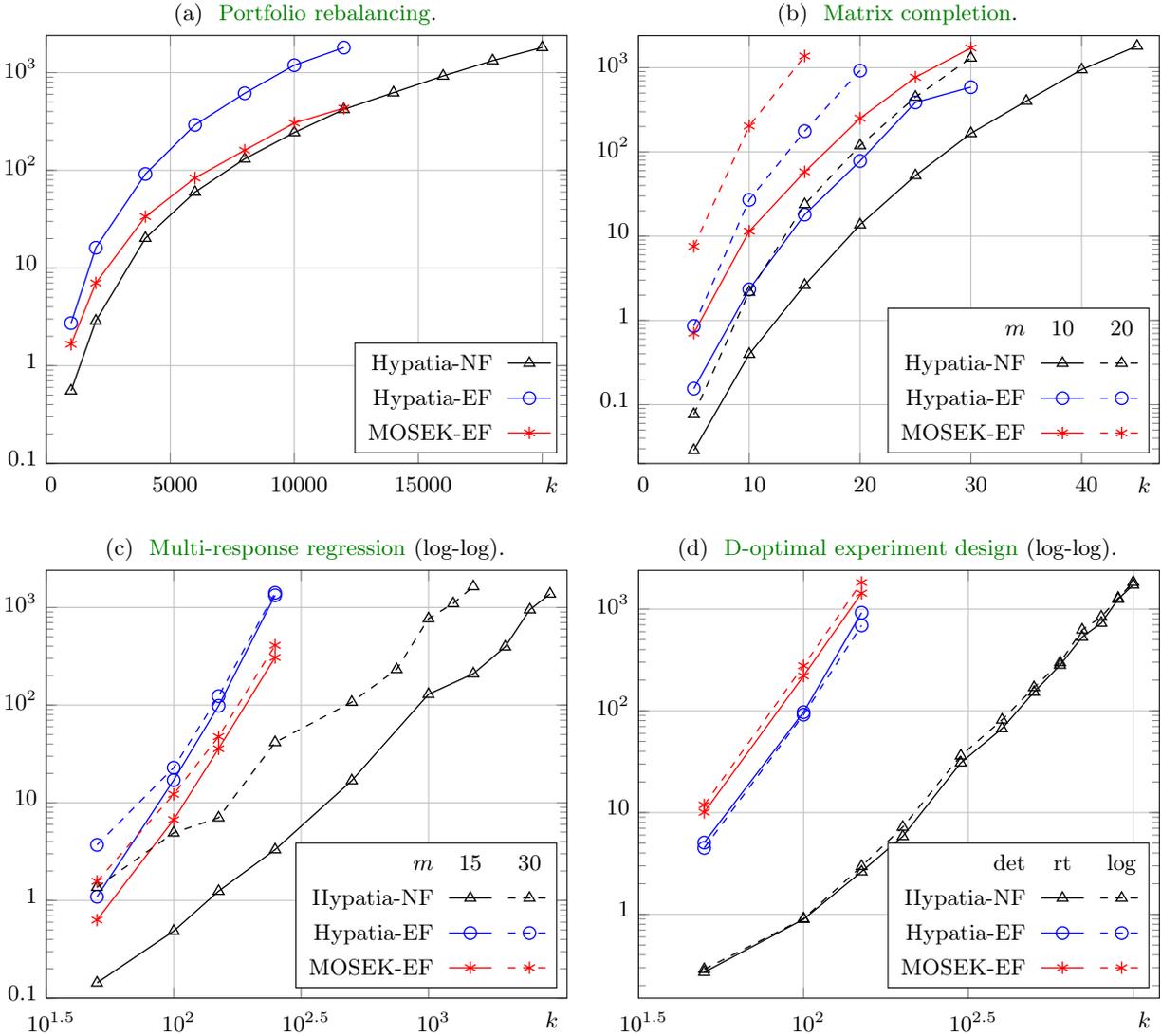
\begin{figure}[!htb]
\centering
\captionsetup{justification=centering}
\caption{
Solve times (in seconds) for solve runs satisfying the convergence check in \cref{eq:convcheck}.
}
\label{fig:plots}
\begin{tikzpicture}[font=\small]
\begin{groupplot}[
    group style={group size=2 by 2, vertical sep=1.5cm, horizontal sep=1cm},
    width = 7.3cm,
    height = 6cm,
    scale only axis,
    every axis plot/.append style={line width=1pt},
    ]
\nextgroupplot[
ymode=log,
ymin=0.1,
ymax=2400,
ytick = {0.1, 1, 10, 100, 1000},
yticklabels={$0.1$, $1$, $10$, $10^2$, $10^3$},
xlabel=$k$,
xmajorgrids,
ymajorgrids,
xmin=0,
xmax=21000,
xtick={0,5000,10000,15000,20000},
xticklabels={$0$,$5000$,$10000$,$15000$},
scaled x ticks=false,
xticklabel style={xshift=+2pt, yshift=-2pt},
yticklabel style={yshift=+2pt},
x label style={at={(axis description cs:1.0,0)}, anchor=north east, yshift=-2pt},
/pgf/number format/.cd, 1000 sep={},
]
\addplot [
black,
line width=0.5pt,
mark=triangle,
mark options={solid},
mark size=2.5pt,
] table [
    x=k,
    y=nat_Hypatia,
    col sep=comma
] {csvs/portfolio_plot.csv};
\label{port:nat}
\addplot [
blue,
line width=0.5pt,
mark=o,
mark options={solid},
mark size=2.5pt,
] table [
    x=k,
    y=SEP_Hypatia,
    col sep=comma
] {csvs/portfolio_plot.csv};
\label{port:ext}
\addplot [
red,
line width=0.5pt,
mark=asterisk,
mark options={solid},
mark size=2.5pt,
] table [
    x=k,
    y=SEP_Mosek,
    col sep=comma
] {csvs/portfolio_plot.csv};
\label{port:mos}
\nextgroupplot[
ymode=log,
ymin=0.02,
ymax=2400,
ytick = {0.01, 0.1, 1, 10, 100, 1000},
yticklabels={, $0.1$, $1$, $10$, $10^2$, $10^3$},
xlabel=$k$,
xmajorgrids,
ymajorgrids,
xmin=0,
xmax=47,
xtick={0,10,20,30,40},
xticklabels={0,10,20,30,40},
xticklabel style={xshift=+2pt, yshift=-2pt},
yticklabel style={yshift=+2pt},
x label style={at={(axis description cs:1.0,0)}, anchor=north east, yshift=-2pt},
]
\addplot [
black,
line width=0.5pt,
mark=triangle,
mark options={solid},
mark size=2.5pt,
] table [
    x=k,
    y=nat_Hypatia,
    col sep=comma
] {csvs/matrixcompletion_plot_10.csv};
\label{matcompl:nat10}
\addplot [
blue,
line width=0.5pt,
mark=o,
mark options={solid},
mark size=2.5pt,
] table [
    x=k,
    y=SEP_Hypatia,
    col sep=comma
] {csvs/matrixcompletion_plot_10.csv};
\label{matcompl:ext10}
\addplot [
red,
line width=0.5pt,
mark=asterisk,
mark options={solid},
mark size=2.5pt,
] table [
    x=k,
    y=SEP_Mosek,
    col sep=comma
] {csvs/matrixcompletion_plot_10.csv};
\label{matcompl:mos10}
\addplot [
dashed,
black,
line width=0.5pt,
mark=triangle,
mark options={solid},
mark size=2.5pt,
] table [
    x=k,
    y=nat_Hypatia,
    col sep=comma
] {csvs/matrixcompletion_plot_20.csv};
\label{matcompl:nat20}
\addplot [
dashed,
blue,
line width=0.5pt,
mark=o,
mark options={solid},
mark size=2.5pt,
] table [
    x=k,
    y=SEP_Hypatia,
    col sep=comma
] {csvs/matrixcompletion_plot_20.csv};
\label{matcompl:ext20}
\addplot [
dashed,
red,
line width=0.5pt,
mark=asterisk,
mark options={solid},
mark size=2.5pt,
] table [
    x=k,
    y=SEP_Mosek,
    col sep=comma
] {csvs/matrixcompletion_plot_20.csv};
\label{matcompl:mos20}
\nextgroupplot[
ymode=log,
ymin=0.1,
ymax=2400,
ytick = {0.1, 1, 10, 100, 1000},
yticklabels={$0.1$,$1$, $10$, $10^2$, $10^3$},
xlabel=$k$,
xmajorgrids,
ymajorgrids,
xmode=log,
xmin=31.6,
xmax=3500,
xtick={31.6, 100, 316.2, 1000},
xticklabels={$10^{1.5}$, $10^2$, $10^{2.5}$, $10^3$},
xticklabel style={xshift=+2pt, yshift=-2pt},
yticklabel style={yshift=+2pt},
x label style={at={(axis description cs:1.0,0)}, anchor=north east, yshift=-2pt},
]
\addplot [
black,
line width=0.5pt,
mark=triangle,
mark options={solid},
mark size=2.5pt,
] table [
    x=k,
    y=nat_Hypatia,
    col sep=comma
] {csvs/matrixregression_plot_15.csv};
\label{matregr:nat15}
\addplot [
blue,
line width=0.5pt,
mark=o,
mark options={solid},
mark size=2.5pt,
] table [
    x=k,
    y=SEP_Hypatia,
    col sep=comma
] {csvs/matrixregression_plot_15.csv};
\label{matregr:ext15}
\addplot [
red,
line width=0.5pt,
mark=asterisk,
mark options={solid},
mark size=2.5pt,
] table [
    x=k,
    y=SEP_Mosek,
    col sep=comma
] {csvs/matrixregression_plot_15.csv};
\label{matregr:mos15}
\addplot [
dashed,
black,
line width=0.5pt,
mark=triangle,
mark options={solid},
mark size=2.5pt,
] table [
    x=k,
    y=nat_Hypatia,
    col sep=comma
] {csvs/matrixregression_plot_30.csv};
\label{matregr:nat30}
\addplot [
dashed,
blue,
line width=0.5pt,
mark=o,
mark options={solid},
mark size=2.5pt,
] table [
    x=k,
    y=SEP_Hypatia,
    col sep=comma
] {csvs/matrixregression_plot_30.csv};
\label{matregr:ext10}
\addplot [
dashed,
red,
line width=0.5pt,
mark=asterisk,
mark options={solid},
mark size=2.5pt,
] table [
    x=k,
    y=SEP_Mosek,
    col sep=comma
] {csvs/matrixregression_plot_30.csv};
\label{matregr:mos30}
\nextgroupplot[
ymode=log,
ymin=0.15,
ymax=2400,
ytick = {0.1, 1, 10, 100, 1000},
yticklabels={$0.1$, $1$, $10$, $10^2$, $10^3$},
xlabel=$k$,
xmajorgrids,
ymajorgrids,
xmode=log,
xmin=31.6,
xmax=1200,
xtick={31.6, 100, 316.2, 1000},
xticklabels={$10^{1.5}$, $10^2$, $10^{2.5}$,},
xticklabel style={xshift=+2pt, yshift=-2pt},
yticklabel style={yshift=+2pt},
x label style={at={(axis description cs:1.0,0)}, anchor=north east, yshift=-2pt},
]
\addplot [
black,
line width=0.5pt,
mark=triangle,
mark options={solid},
mark size=2.5pt,
] table [
    x=k,
    y=nat_Hypatia,
    col sep=comma
] {csvs/doptimaldesign_plot_false.csv};
\label{expdesign:natrt}
\addplot [
blue,
line width=0.5pt,
mark=o,
mark options={solid},
mark size=2.5pt,
] table [
    x=k,
    y=SEP_Hypatia,
    col sep=comma
] {csvs/doptimaldesign_plot_false.csv};
\label{expdesign:extrt}
\addplot [
red,
line width=0.5pt,
mark=asterisk,
mark options={solid},
mark size=2.5pt,
] table [
    x=k,
    y=SEP_Mosek,
    col sep=comma
] {csvs/doptimaldesign_plot_false.csv};
\label{expdesign:mosrt}
\addplot [
dashed,
black,
line width=0.5pt,
mark=triangle,
mark options={solid},
mark size=2.5pt,
] table [
    x=k,
    y=nat_Hypatia,
    col sep=comma
] {csvs/doptimaldesign_plot_true.csv};
\label{expdesign:natlog}
\addplot [
dashed,
blue,
line width=0.5pt,
mark=o,
mark options={solid},
mark size=2.5pt,
] table [
    x=k,
    y=EP_Hypatia,
    col sep=comma
] {csvs/doptimaldesign_plot_true.csv};
\label{expdesign:extlog}
\addplot [
dashed,
red,
line width=0.5pt,
mark=asterisk,
mark options={solid},
mark size=2.5pt,
] table [
    x=k,
    y=EP_Mosek,
    col sep=comma
] {csvs/doptimaldesign_plot_true.csv};
\label{expdesign:moslog}
\end{groupplot}
\node [text width=15em, above=-0.1em] at (group c1r1.north)
{\subcaption{\label{fig:plots:portfolio} \nameref{sec:testing:portfolio}.}};
\node [text width=15em, above=-0.1em] at (group c2r1.north)
{\subcaption{\label{fig:plots:matrixcompletion} \nameref{sec:testing:matrixcompletion}.}};
\node [text width=20em, above=-0.1em] at (group c1r2.north)
{\subcaption{\label{fig:plots:matrixregression} \nameref{sec:testing:matrixregression} (log-log).}};
\node [text width=20em, above=-0.1em] at (group c2r2.north)
{\subcaption{\label{fig:plots:expdesign} \nameref{sec:testing:expdesign} (log-log).}};
\node[
draw,
fill=white,
inner sep=1pt,
inner xsep=-1pt,
above left=0.4em,
] at (group c1r1.south east) {
\setlength{\tabcolsep}{3pt}
\renewcommand{\arraystretch}{1.1}
\begin{tabular}{@{}rc}
Hypatia-NF & \ref{port:nat} \\
Hypatia-EF & \ref{port:ext} \\
MOSEK-EF & \ref{port:mos}
\end{tabular}
};
\node[
draw,
fill=white,
inner sep=1pt,
inner xsep=-1pt,
above left=0.4em,
] at (group c2r1.south east) {
\setlength{\tabcolsep}{3pt}
\renewcommand{\arraystretch}{1.1}
\begin{tabular}{@{}rcc}
$m$ & 10 & 20 \\
Hypatia-NF & \ref{matcompl:nat10} & \ref{matcompl:nat20} \\
Hypatia-EF & \ref{matcompl:ext10} & \ref{matcompl:ext20} \\
MOSEK-EF & \ref{matcompl:mos10} & \ref{matcompl:mos20}
\end{tabular}
};
\node[
draw,
fill=white,
inner sep=1pt,
inner xsep=-1pt,
above left=0.4em,
] at (group c1r2.south east) {
\setlength{\tabcolsep}{3pt}
\renewcommand{\arraystretch}{1.1}
\begin{tabular}{@{}rcc}
$m$ & 15 & 30 \\ 
Hypatia-NF & \ref{matregr:nat15} & \ref{matregr:nat30} \\
Hypatia-EF & \ref{matregr:ext15} & \ref{matregr:ext10} \\
MOSEK-EF & \ref{matregr:mos15} & \ref{matregr:mos30}
\end{tabular}
};
\node[
draw,
fill=white,
inner sep=1pt,
inner xsep=-1pt,
above left=0.4em,
] at (group c2r2.south east) {
\setlength{\tabcolsep}{3pt}
\renewcommand{\arraystretch}{1.1}
\begin{tabular}{@{}rcc}
det & rt & log \\
Hypatia-NF & \ref{expdesign:natrt} & \ref{expdesign:natlog} \\
Hypatia-EF & \ref{expdesign:extrt} & \ref{expdesign:extlog} \\
MOSEK-EF & \ref{expdesign:mosrt} & \ref{expdesign:moslog}
\end{tabular}
};
\end{tikzpicture}

\end{figure}

In \crefrange{sec:testing:portfolio}{sec:testing:shapeconregr}, we present example problems with NFs using some of Hypatia's predefined cones and EFs constructed using the techniques from \cref{sec:cones}. 
For each example problem, we generate random instances of a wide variety of sizes, and we observe larger dimensions and often larger barrier parameters for EFs compared to NFs.
In \crefrange{tab:portfolio}{tab:shapeconregr}, $\nu$ and $n$, $p$, $q$ refer to the NF barrier parameter and primal variable, linear equality, and cone inequality dimensions (in our general conic form \cref{eq:prim}), and $\bar{\nu}$, $\bar{n}$, $\bar{p}$, $\bar{q}$ refer to the corresponding EF values.
For three solver/formulation combinations - Hypatia with NF (\emph{Hypatia-NF}), Hypatia with EF (\emph{Hypatia-EF}), and MOSEK with EF (\emph{MOSEK-EF}) - we compare termination statuses, iteration counts, and solve times in seconds (columns \emph{st}, \emph{it}, and \emph{time}) in \crefrange{tab:portfolio}{tab:shapeconregr} and \cref{fig:plots}.
In \cref{sec:testing:matrixcompletion,sec:testing:expdesign} we depend on a geometric mean cone EF, so we compare the \emph{EF-exp} and \emph{EF-sec} formulations from \cref{sec:cones:geom}.
Note that all of our instances are primal-dual feasible, so we expect solvers to return optimality certificates.
Compared to Hypatia-EF and MOSEK-EF, Hypatia-NF generally converges faster and more reliably, and solves larger instances within time and memory limits.

We perform all instance generation, computational experiments, and results analysis with Ubuntu 21.04, Julia 1.7, and Hypatia 0.5.0, on dedicated hardware with an AMD Ryzen 9 3950X 16-core processor (32 threads) and 128GB of RAM.
We limit each solver to using 16 threads.
We use JuMP 0.21.5 and MathOptInterface 0.9.18 to build all instances.
We use MOSEK 9 through MosekTools.jl 0.9.4 (which is maintained in part by MOSEK).\footnote{
We note that MOSEK 9's primal conic form only recognizes conic constraints of the form $x \in \K$ \citep[Section 8]{mosek2020modeling}, whereas Hypatia accepts the more general affine form $h - Gx \in \K$ (see \cref{eq:prim:K}).
MathOptInterface recognizes both \emph{VectorOfVariables} form $x \in \K$ and \emph{VectorAffineFunction} form $h - Gx \in \K$.
Since JuMP and MathOptInterface (including bridges) use the $x \in \K$ form whenever possible, unnecessary high dimensional slack variables are not introduced when instances in Hypatia's general conic form \cref{eq:prim} are converted into MOSEK 9's form.
}
MOSEK uses its conic interior point method for all solves.
We note that MOSEK heuristically determines whether it is more efficient to solve the primal or dual of an instance \citep[Section 13.1]{aps2020mosek}; Hypatia does not do this.
We do not disable any MOSEK features.
Hypatia uses one particular default algorithmic implementation that we describe broadly in \citet{coey2021performance} (the \emph{combined directions method}); this is not the most efficient method for all instances, but Hypatia does not currently have heuristics for choosing which stepper or linear system solver procedure to use, for example.
Simple scripts and instructions for reproducing all results are available in Hypatia's benchmarks/natvsext folder.
A CSV file containing raw results is available at the Hypatia wiki page. 

Hypatia and MOSEK use similar convergence criteria (see \citet[Section 13.3.2]{aps2020mosek}), and we set their feasibility and optimality gap tolerances to $10^{-7}$.
In the solver statistics tables, asterisks indicate missing data, and we use the following codes for the termination status (st) columns:
\begin{description}[font=\normalfont\itshape]
\item[co]- the solver claims it has an approximate optimality certificate,
\item[tl]- the solver stops itself due to a solve time limit of 1800 seconds, or the solve run is killed because it takes at least $1.2 \times 1800$ seconds,
\item[rl]- the solve is terminated because insufficient RAM is available,
\item[sp]- the solver reports \emph{slow progress} during iterations,
\item[er]- the solver encounters a numerical error,
\item[m]- the model cannot be constructed with JuMP due to insufficient RAM or a model generation time limit of $1.2 \times 1800$ seconds (EFs tend to be slower and more memory-intensive to construct than NFs, so EF columns often have missing data),
\item[sk]- we skip the solve run because a smaller instance has a tl or rl status, or we skip model generation because a smaller instance has an m status.
\end{description}

For each solve run that yields a primal-dual point $(x, y, z, s)$ (see \cref{sec:form}; $s \in \K$ and $z \in \K^\ast$ are the solver's primal and dual cone interior points at termination), we compute:
\begin{equation}
\epsilon \coloneqq \max \biggl\{ 
\frac{\lVert A' y + G' z + c \rVert_{\infty}}{1 + \norm{c}_{\infty}}, 
\frac{\norm{-A x + b}_{\infty}}{1 + \norm{b}_{\infty}}, 
\frac{\norm{-G x + h - s}_{\infty}}{1 + \norm{h}_{\infty}},
\frac{\lvert c' x + b' y + h' z \rvert}{1 + \lvert b' y + h' z \rvert}
\biggr\},
\label{eq:convcheck}
\end{equation}
and if $\epsilon < 10^{-5}$, we underline the corresponding status code (e.g. \underline{co}, \underline{tl}) to indicate that the solution approximately satisfies the optimality certificate conditions from \cref{sec:form}.
In our solve time plots in \cref{fig:plots}, we only plot solve runs with underlined status codes.
Finally, for each instance and each pair of corresponding solve runs with \underline{co} status codes, we compute the relative difference of the primal objective values $g_1$ and $g_2$ as $\tilde{\epsilon} \coloneqq \lvert g_1 - g_2 \rvert / (1 + \max(\lvert g_1 \rvert, \lvert g_2 \rvert))$.
We note $\tilde{\epsilon} < 10^{-5}$ for most instances and pairs of solvers, and $\tilde{\epsilon} < 10^{-3}$ in all cases.

\subsection{Portfolio rebalancing}
\label{sec:testing:portfolio}

Suppose there are $k$ possible investments with expected returns $g \in \bbR_>^k$ and covariance matrix $\Sigma \in \bbS_{\succ}^k$.
We let $\rho \in [-1, 1]^k$ be the investment variable, which must also satisfy side constraints $F \rho = 0$, where $F \in \bbR^{l \times k}$.
We formulate a risk-constrained portfolio rebalancing optimization problem as:
\begin{subequations}
\begin{align}
\textstyle\max_{\rho \in \bbR^k} \quad g' \rho & :
\\
e' \rho & = 0,
\\
F \rho & = 0,
\\
( 1, \rho ) & \in \klinf,
\label{eq:portfolio:inf}
\\
( \gamma, \Sigma^{1/2} \rho ) & \in \klinf^\ast.
\label{eq:portfolio:1} 
\end{align}
\label{eq:portfolio}
\end{subequations}
Note \cref{eq:portfolio:inf} expresses $\rho \in [-1, 1]^k$ and \cref{eq:portfolio:1} is a risk constraint.
The EFs for \cref{eq:portfolio:inf,eq:portfolio:1} follow \cref{eq:cones:linf,eq:cones:l1}. 
Note the EF is a standard linear program.

To build random instances of \cref{eq:portfolio}, we generate $g$ with independent uniform positive entries, and $\Sigma^{1/2}$ and $F$ with independent Gaussian entries, for $l = k / 2$ and various values of $k$.
We use $\Sigma^{1/2}$ to compute reasonable values for the risk parameter $\gamma > 0$. 
Our results are summarized in \cref{tab:portfolio,fig:plots:portfolio}.
Note that $\nu = q = 2k + 2$, $\bar{\nu} = \bar{q} = 4k + 1$, $n = k$, $\bar{n} = 2k$, $p = \bar{p} = k/2 + 1$.
The variable and conic constraint dimensions of the EFs are approximately double those of the NFs. 
Hypatia-NF exhibits the fastest solve times and solves much larger instances than Hypatia-EF and MOSEK-EF.
MOSEK requires notably fewer PDIPM iterations than Hypatia.

\begin{table}[!htb]
\centering
\caption{
\nameref{sec:testing:portfolio} solver statistics.
}
\label{tab:portfolio}
\small
\sisetup{
table-text-alignment = right,
table-auto-round,
table-figures-integer = 4,
table-figures-decimal = 1,
table-format = 4.1,
add-decimal-zero = false,
add-integer-zero = false,
}
\begin{tabular}{rcrScrScrS}
\toprule
& \multicolumn{3}{c}{Hypatia-NF} & \multicolumn{3}{c}{Hypatia-EF} & \multicolumn{3}{c}{MOSEK-EF} \\
\cmidrule(lr){2-4} \cmidrule(lr){5-7} \cmidrule(lr){8-10} 
$k$ & {st} & {it} & {time} & {st} & {it} & {time} & {st} & {it} & {time} \\
\midrule
1000 & \underline{co} & 31 & 0.6 & \underline{co} & 25 & 2.7 & \underline{co} & 9 & 1.7 \\
2000 & \underline{co} & 36 & 2.9 & \underline{co} & 28 & 16. & \underline{co} & 10 & 7.0 \\
4000 & \underline{co} & 45 & 20. & \underline{co} & 29 & 92. & \underline{co} & 10 & 34. \\
6000 & \underline{co} & 49 & 60. & \underline{co} & 34 & 292. & \underline{co} & 10 & 83. \\
8000 & \underline{co} & 51 & 131. & \underline{co} & 33 & 615. & \underline{co} & 10 & 160. \\
10000 & \underline{co} & 55 & 244. & \underline{co} & 36 & 1192. & \underline{co} & 12 & 305. \\
12000 & \underline{co} & 62 & 421. & \underline{tl} & 32 & 1805. & \underline{co} & 10 & 433. \\
14000 & \underline{co} & 61 & 624. & sk & $\ast$ & $\ast$ & rl & $\ast$ & $\ast$ \\
16000 & \underline{co} & 63 & 924. & sk & $\ast$ & $\ast$ & sk & $\ast$ & $\ast$ \\
18000 & \underline{co} & 64 & 1327. & sk & $\ast$ & $\ast$ & sk & $\ast$ & $\ast$ \\
20000 & \underline{co} & 66 & 1810. & sk & $\ast$ & $\ast$ & sk & $\ast$ & $\ast$ \\
\bottomrule
\end{tabular}
\end{table}

\subsection{Matrix completion}
\label{sec:testing:matrixcompletion}

Suppose there exists a matrix $F \in \bbR^{k \times l}$ and we know the entries $(F_{i,j})_{(i,j) \in \mathcal{S}}$ in the sparsity pattern $\mathcal{S}$.
In the matrix completion problem, we seek to estimate the missing components $(F_{i,j})_{(i,j) \not\in \mathcal{S}}$. 
We modify the formulation in \citet[Section 4.3]{agrawal2019disciplined} by replacing the spectral radius in the objective function with the spectral norm (allowing rectangular matrices) and using a convex relaxation of the geometric mean equality constraint:
\begin{subequations}
\begin{align}
\textstyle\min_{\rho \in \bbR, X \in \bbR^{k \times l}} \quad \rho &:
\\
X_{i,j} &= F_{i,j} 
\quad \forall (i,j) \in \mathcal{S},
\\
(\rho, \vect(X)) & \in \K_{\lspec(k, l)},
\label{eq:matrixcompletion:spec}
\\
\bigl( 1, (X_{i,j})_{(i,j) \not\in \mathcal{S}} \bigr) & \in \kgeom.
\label{eq:matrixcompletion:geom}
\end{align}
\label{eq:matrixcompletion}
\end{subequations}
The EF for \cref{eq:matrixcompletion:spec} follows \cref{eq:cones:spec}, and for \cref{eq:matrixcompletion:geom} we compare EF-exp and EF-sec (see \cref{sec:cones:geom}).

To build random instances of \cref{eq:matrixcompletion}, we generate sparse matrices $F$ with independent Gaussian nonzero entries, for various values of $k$, column-to-row ratios $m \in \{10, 20\}$, and $l = m k$.
Our results are summarized in \cref{tab:matrixcompletion:sizes,tab:matrixcompletion:times,fig:plots:matrixcompletion}.
Note we only plot EF-sec results for Hypatia-EF and MOSEK-EF, as MOSEK performs better with EF-sec (which only uses symmetric cones) than with EF-exp, though Hypatia exhibits the opposite trend.
Hypatia-NF is much faster and solves more instances than the Hypatia-EFs and MOSEK-EFs.

\begin{table}[!htb]
\centering
\caption{
\nameref{sec:testing:matrixcompletion} formulation statistics. 
Note $p = \bar{p} = \lvert S \rvert$.
}
\label{tab:matrixcompletion:sizes}
\small
\begin{tabular}{rrrrrrrrrrrr}
\toprule
& & \multicolumn{4}{c}{NF} & \multicolumn{3}{c}{EF-exp} & \multicolumn{3}{c}{EF-sec} \\
\cmidrule(lr){3-6} \cmidrule(lr){7-9} \cmidrule(lr){10-12}
$m$ & $k$ & $\nu$ & $n$ & $p$ & $q$ & $\bar{\nu}$ & $\bar{n}$ & $\bar{q}$ & $\bar{\nu}$ & $\bar{n}$ & $\bar{q}$ \\
\midrule
\multirow{9}{*}{10}
& 5 & 57 & 251 & 200 & 302 & 207 & 302 & 1692 & 182 & 314 & 1730 \\
& 10 & 218 & 1001 & 794 & 1208 & 730 & 1208 & 6725 & 621 & 1256 & 6871 \\
& 15 & 472 & 2251 & 1795 & 2707 & 1532 & 2707 & 15062 & 1188 & 2762 & 15229 \\
& 20 & 846 & 4001 & 3176 & 4826 & 2694 & 4826 & 26784 & 2267 & 5024 & 27380 \\
& 25 & 1299 & 6251 & 4978 & 7524 & 4093 & 7524 & 41768 & 4370 & 8298 & 44092 \\
& 30 & 1858 & 9001 & 7174 & 10828 & 5810 & 10828 & 60095 & 4425 & 11048 & 60757 \\
& 35 & 2477 & 12251 & 9810 & 14692 & 7707 & 14692 & 81627 & 8576 & 16346 & 86591 \\
& 40 & 3256 & 16001 & 12786 & 19216 & 10084 & 19216 & 106664 & 8631 & 20096 & 109306 \\
& 45 & 4142 & 20251 & 16155 & 24347 & 12782 & 24347 & 135047 & 8686 & 24346 & 135046 \\
\addlinespace[2pt]
\hdashline
\addlinespace[2pt]
\multirow{6}{*}{20}
& 5 & 114 & 501 & 393 & 609 & 428 & 609 & 5888 & 360 & 628 & 5947 \\
& 10 & 418 & 2001 & 1594 & 2408 & 1430 & 2408 & 23375 & 1233 & 2512 & 23689 \\
& 15 & 933 & 4501 & 3584 & 5418 & 3065 & 5418 & 52520 & 2362 & 5524 & 52840 \\
& 20 & 1663 & 8001 & 6359 & 9643 & 5345 & 9643 & 93335 & 4515 & 10048 & 94552 \\
& 25 & 2513 & 12501 & 10014 & 14988 & 7985 & 14988 & 145535 & 8716 & 16596 & 150361 \\
& 30 & 3643 & 18001 & 14389 & 21613 & 11465 & 21613 & 209600 & 8821 & 22096 & 211051 \\
\bottomrule
\end{tabular}
\end{table}

\begin{table}[!htb]
\centering
\caption{
\nameref{sec:testing:matrixcompletion} solver statistics.
}
\label{tab:matrixcompletion:times}
\small
\sisetup{
table-text-alignment = right,
table-auto-round,
table-figures-integer = 4,
table-figures-decimal = 1,
table-format = 4.1,
add-decimal-zero = false,
add-integer-zero = false,
}
\begin{tabular}{rrcrScrScrScrScrS}
\toprule
& & \multicolumn{3}{c}{NF} & \multicolumn{6}{c}{EF-exp} & \multicolumn{6}{c}{EF-sec} \\
\cmidrule(lr){3-5} \cmidrule(lr){6-11} \cmidrule(lr){12-17}
& & \multicolumn{3}{c}{Hypatia} & \multicolumn{3}{c}{Hypatia} & \multicolumn{3}{c}{MOSEK} & \multicolumn{3}{c}{Hypatia} & \multicolumn{3}{c}{MOSEK} \\
\cmidrule(lr){3-5} \cmidrule(lr){6-8} \cmidrule(lr){9-11} \cmidrule(lr){12-14} \cmidrule(lr){15-17} 
$m$ & $k$ & {st} & {it} & {time} & {st} & {it} & {time} & {st} & {it} & {time} & {st} & {it} & {time} & {st} & {it} & {time} \\
\midrule
\multirow{9}{*}{10}
& 5 & \underline{co} & 14 & 0.0 & \underline{co} & 19 & 0.1 & \underline{co} & 15 & 0.9 & \underline{co} & 18 & 0.2 & \underline{co} & 11 & 0.7 \\
& 10 & \underline{co} & 19 & 0.4 & \underline{co} & 34 & 2.2 & \underline{co} & 20 & 20. & \underline{co} & 30 & 2.3 & \underline{co} & 10 & 11. \\
& 15 & \underline{co} & 23 & 2.6 & \underline{co} & 42 & 16. & \underline{co} & 21 & 120. & \underline{co} & 41 & 18. & \underline{co} & 9 & 58. \\
& 20 & \underline{co} & 26 & 14. & \underline{co} & 52 & 70. & \underline{co} & 24 & 524. & \underline{co} & 47 & 78. & \underline{co} & 11 & 251. \\
& 25 & \underline{co} & 30 & 52. & \underline{co} & 59 & 225. & \underline{co} & 26 & 1758. & \underline{co} & 69 & 387. & \underline{co} & 11 & 770. \\
& 30 & \underline{co} & 34 & 166. & \underline{co} & 61 & 556. & tl & 11 & 1817. & \underline{co} & 55 & 587. & \underline{co} & 10 & 1712. \\
& 35 & \underline{co} & 39 & 402. & \underline{co} & 61 & 1228. & sk & $\ast$ & $\ast$ & tl & 61 & 1817. & rl & $\ast$ & $\ast$ \\
& 40 & \underline{co} & 48 & 949. & tl & 34 & 1820. & sk & $\ast$ & $\ast$ & sk & $\ast$ & $\ast$ & sk & $\ast$ & $\ast$ \\
& 45 & \underline{co} & 47 & 1806. & sk & $\ast$ & $\ast$ & sk & $\ast$ & $\ast$ & sk & $\ast$ & $\ast$ & sk & $\ast$ & $\ast$ \\
\addlinespace[2pt]
\hdashline
\addlinespace[2pt]
\multirow{6}{*}{20}
& 5 & \underline{co} & 15 & 0.1 & \underline{co} & 29 & 0.8 & \underline{co} & 17 & 14. & \underline{co} & 26 & 0.9 & \underline{co} & 8 & 7.6 \\
& 10 & \underline{co} & 22 & 2.2 & \underline{co} & 48 & 25. & \underline{co} & 25 & 448. & \underline{co} & 45 & 27. & \underline{co} & 10 & 203. \\
& 15 & \underline{co} & 30 & 24. & \underline{co} & 59 & 179. & tl & 14 & 1871. & \underline{co} & 52 & 176. & \underline{co} & 10 & 1375. \\
& 20 & \underline{co} & 33 & 119. & \underline{co} & 71 & 786. & sk & $\ast$ & $\ast$ & \underline{co} & 70 & 924. & rl & $\ast$ & $\ast$ \\
& 25 & \underline{co} & 41 & 448. & tl & 47 & 1822. & sk & $\ast$ & $\ast$ & tl & 26 & 1804. & sk & $\ast$ & $\ast$ \\
& 30 & \underline{co} & 52 & 1305. & sk & $\ast$ & $\ast$ & sk & $\ast$ & $\ast$ & sk & $\ast$ & $\ast$ & sk & $\ast$ & $\ast$ \\
\bottomrule
\end{tabular}
\end{table}

\subsection{Multi-response regression}
\label{sec:testing:matrixregression}

In the multi-response linear regression problem, we seek to estimate a coefficient matrix $F \in \bbR^{m \times l}$ from a design matrix $X \in \bbR^{l \times k}$ and response matrix $Y \in \bbR^{m \times k}$.
We use a similar formulation to the one proposed in \citet{yang2016nuclear}, with nuclear norm loss and $\ell_2$ norm regularization:
\begin{subequations}
\begin{align}
\textstyle\min_{\rho \in \bbR, \mu \in \bbR, F \in \bbR^{m \times l}} \quad \rho + \gamma \mu & :
\\
(\rho, \vect(Y - F X)) & \in \K_{\lspec(m, k)}^\ast,
\label{eq:matrixregression:nuc}
\\
(\mu, \vect(F)) & \in \kltwo.
\end{align}
\label{eq:matrixregression}
\end{subequations}
The EF for NF constraint \cref{eq:matrixregression:nuc} follows \cref{eq:cones:nuc}. 

To build random instances of \cref{eq:matrixregression}, we generate random $X$ and $Y$ with independent Gaussian entries, for various values of $k$ with $l = m \in \{15, 30\}$, and we use regularization parameter $\gamma = 0.1$.
Our results are summarized in \cref{tab:matrixregression,fig:plots:matrixregression}.
Note that $\nu = 3 + m$, $\bar{\nu} = \nu + k$, $n = 2 + m^2$, $p = \bar{p} = 0$, $\bar{q} = \bar{n} + m k$.
The variable dimensions for the NFs only depend on $k$ and are much smaller than those of the EFs. 
The EFs also have much larger conic constraint dimensions.
Hypatia-NF exhibits faster solve times than Hypatia-EF and MOSEK-EF. 
Hypatia-NF solves much larger instances and takes a fairly consistent number of iterations.

\begin{table}[!htb]
\centering
\caption{
\nameref{sec:testing:matrixregression} formulation and solver statistics. 
}
\label{tab:matrixregression}
\small
\sisetup{
table-text-alignment = right,
table-auto-round,
table-figures-integer = 4,
table-figures-decimal = 1,
table-format = 4.1,
add-decimal-zero = false,
add-integer-zero = false,
}
\begin{tabular}{rrrrcrScrScrS}
\toprule
& & \multicolumn{2}{c}{form. stats.} & \multicolumn{3}{c}{Hypatia-NF} & \multicolumn{3}{c}{Hypatia-EF} & \multicolumn{3}{c}{MOSEK-EF} \\
\cmidrule(lr){3-4} \cmidrule(lr){5-7} \cmidrule(lr){8-10} \cmidrule(lr){11-13} 
$m$ & $k$ & $\bar{n}$ & $q$ & {st} & {it} & {time} & {st} & {it} & {time} & {st} & {it} & {time} \\
\midrule
\multirow{10}{*}{15}
& 50 & 1622 & 977 & \underline{co} & 11 & 0.1 & \underline{co} & 12 & 1.1 & \underline{co} & 4 & 0.6 \\
& 100 & 5397 & 1727 & \underline{co} & 10 & 0.5 & \underline{co} & 12 & 17. & \underline{co} & 5 & 6.7 \\
& 150 & 11672 & 2477 & \underline{co} & 10 & 1.2 & \underline{co} & 13 & 98. & \underline{co} & 5 & 36. \\
& 250 & 31722 & 3977 & \underline{co} & 10 & 3.3 & \underline{co} & 14 & 1331. & \underline{co} & 5 & 308. \\
& 500 & 125597 & 7727 & \underline{co} & 10 & 17. & m & $\ast$ & $\ast$ & tl & $\ast$ & $\ast$ \\
& 1000 & $\ast$ & 15227 & \underline{co} & 13 & 129. & sk & $\ast$ & $\ast$ & sk & $\ast$ & $\ast$ \\
& 1500 & $\ast$ & 22727 & \underline{co} & 10 & 209. & sk & $\ast$ & $\ast$ & sk & $\ast$ & $\ast$ \\
& 2000 & $\ast$ & 30227 & \underline{co} & 9 & 395. & sk & $\ast$ & $\ast$ & sk & $\ast$ & $\ast$ \\
& 2500 & $\ast$ & 37727 & \underline{co} & 11 & 949. & sk & $\ast$ & $\ast$ & sk & $\ast$ & $\ast$ \\
& 3000 & $\ast$ & 45227 & \underline{co} & 11 & 1375. & sk & $\ast$ & $\ast$ & sk & $\ast$ & $\ast$ \\
\addlinespace[2pt]
\hdashline
\addlinespace[2pt]
\multirow{9}{*}{30}
& 50 & 2642 & 2402 & \underline{co} & 13 & 1.4 & \underline{co} & 11 & 3.7 & \underline{co} & 5 & 1.6 \\
& 100 & 6417 & 3902 & \underline{co} & 14 & 4.9 & \underline{co} & 11 & 23. & \underline{co} & 5 & 12. \\
& 150 & 12692 & 5402 & \underline{co} & 12 & 7.0 & \underline{co} & 12 & 123. & \underline{co} & 5 & 47. \\
& 250 & 32742 & 8402 & \underline{co} & 15 & 41. & \underline{co} & 13 & 1412. & \underline{co} & 5 & 409. \\
& 500 & 126617 & 15902 & \underline{co} & 11 & 107. & m & $\ast$ & $\ast$ & tl & $\ast$ & $\ast$ \\
& 750 & $\ast$ & 23402 & \underline{co} & 11 & 232. & sk & $\ast$ & $\ast$ & sk & $\ast$ & $\ast$ \\
& 1000 & $\ast$ & 30902 & \underline{co} & 13 & 768. & sk & $\ast$ & $\ast$ & sk & $\ast$ & $\ast$ \\
& 1250 & $\ast$ & 38402 & \underline{co} & 12 & 1098. & sk & $\ast$ & $\ast$ & sk & $\ast$ & $\ast$ \\
& 1500 & $\ast$ & 45902 & \underline{co} & 12 & 1637. & sk & $\ast$ & $\ast$ & sk & $\ast$ & $\ast$ \\
\bottomrule
\end{tabular}
\end{table}

\subsection{D-optimal experiment design}
\label{sec:testing:expdesign}

In a continuous relaxation of the D-optimal experiment design problem (see \citet[Section 7.5]{boyd2004convex}), the variable $\mu \in \bbR^m$ is the number of trials to run for each of $m$ experiments, and our goal is to minimize the determinant of the error covariance matrix $(F \Diag(\mu) F')^{-1}$, given a menu of experiments $F \in \bbR^{k \times m}$ useful for estimating a vector in $\bbR^k$.
We require that a total of $j$ experiments are performed and that each experiment can be performed between $0$ and $l$ times.
We formulate this problem as:
\begin{subequations}
\begin{align}
\textstyle\max_{\rho \in \bbR, \mu \in \bbR^m} \quad \rho &:
\\
e' \mu & = j,
\\
( l / 2, \mu - (l / 2) e ) & \in \klinf,
\label{eq:expdesign:infty}
\\
( \rho, \vect ( F \Diag(\mu) F' ) ) & \in \krtdet.
\label{eq:expdesign:rootdet}
\end{align}
\label{eq:expdesign}
\end{subequations}
In an alternative \emph{logdet} variant of the \emph{rtdet} variant \cref{eq:expdesign}, we replace \cref{eq:expdesign:rootdet} with:
\begin{equation}
( \rho, 1, \vect ( F \Diag(\mu) F' ) ) \in \klogdet,
\label{eq:expdesign:logdet}
\end{equation}
noting that both variants have the same optimal solution set for $\mu$.
The EFs for \cref{eq:expdesign:infty,eq:expdesign:rootdet,eq:expdesign:logdet} follow \cref{eq:cones:linf,eq:cones:rootdet,eq:cones:logdet}.
Since the EF for $\krtdet$ depends on a $\kgeom$ EF, for the rtdet variant we compare EF-exp and EF-sec (see \cref{sec:cones:geom}).

To build random instances of \cref{eq:expdesign}, we generate $F$ with independent Gaussian entries, for various values of $k$, $m = j = 2k$, and $l = 5$.
Our results are summarized in \cref{tab:expdesign:logdet,tab:expdesign:rtdettimes,fig:plots:expdesign}.
For the logdet variant, $\nu = 3 + 3k$, $\bar{\nu} = 1 + 9k$, $n = 1 + 2k$, $p = \bar{p} = 1$.
The sizes for the rootdet formulations are similar to those of the logdet formulations, so we exclude these.
Note for the rootdet variant, we only plot EF-sec results for Hypatia-EF and MOSEK-EF, as MOSEK typically performs slightly better with EF-sec than with EF-exp, though Hypatia exhibits the opposite trend.
For both variants, the NFs have much lower variable and conic constraint dimensions than the EFs.
Although the EF solvers typically solve instances up to $k = 150$, Hypatia-NF solves instances with $k$ at least $900$.
Hypatia-NF is also much faster than the EF solvers for all $k$.

\begin{table}[!htb]
\centering
\caption{
\nameref{sec:testing:expdesign} logdet variant formulation and solver statistics. 
}
\label{tab:expdesign:logdet}
\small
\sisetup{
table-text-alignment = right,
table-auto-round,
table-figures-integer = 4,
table-figures-decimal = 1,
table-format = 4.1,
add-decimal-zero = false,
add-integer-zero = false,
}
\begin{tabular}{rrrrcrScrScrS}
\toprule
& \multicolumn{3}{c}{form. stats.} & \multicolumn{3}{c}{Hypatia-NF} & \multicolumn{3}{c}{Hypatia-EF} & \multicolumn{3}{c}{MOSEK-EF} \\
\cmidrule(lr){2-4} \cmidrule(lr){5-7} \cmidrule(lr){8-10} \cmidrule(lr){11-13} 
$k$ & $\bar{n}$ & $q$ & $\bar{q}$ & {st} & {it} & {time} & {st} & {it} & {time} & {st} & {it} & {time} \\
\midrule
50 & 1426 & 1378 & 5401 & \underline{co} & 25 & 0.3 & \underline{co} & 21 & 4.5 & \underline{co} & 15 & 12. \\
100 & 5351 & 5253 & 20801 & \underline{co} & 26 & 0.9 & \underline{co} & 25 & 91. & \underline{co} & 15 & 277. \\
150 & 11776 & 11628 & 46201 & \underline{co} & 29 & 3.0 & \underline{co} & 27 & 690. & \underline{tl} & 14 & 1825. \\
200 & 20701 & 20503 & 81601 & \underline{co} & 28 & 7.2 & tl & 17 & 1849. & sk & $\ast$ & $\ast$ \\
300 & 46051 & 45753 & 182401 & \underline{co} & 36 & 36. & sk & $\ast$ & $\ast$ & sk & $\ast$ & $\ast$ \\
400 & 81401 & 81003 & 323201 & \underline{co} & 36 & 81. & m & $\ast$ & $\ast$ & sk & $\ast$ & $\ast$ \\
500 & 126751 & 126253 & 504001 & \underline{co} & 36 & 169. & sk & $\ast$ & $\ast$ & sk & $\ast$ & $\ast$ \\
600 & 182101 & 181503 & 724801 & \underline{co} & 36 & 298. & sk & $\ast$ & $\ast$ & sk & $\ast$ & $\ast$ \\
700 & $\ast$ & 246753 & $\ast$ & \underline{co} & 39 & 624. & sk & $\ast$ & $\ast$ & m & $\ast$ & $\ast$ \\
800 & $\ast$ & 322003 & $\ast$ & \underline{co} & 37 & 838. & sk & $\ast$ & $\ast$ & sk & $\ast$ & $\ast$ \\
900 & $\ast$ & 407253 & $\ast$ & \underline{co} & 37 & 1282. & sk & $\ast$ & $\ast$ & sk & $\ast$ & $\ast$ \\
1000 & $\ast$ & 502503 & $\ast$ & \underline{tl} & 37 & 1838. & sk & $\ast$ & $\ast$ & sk & $\ast$ & $\ast$ \\
\bottomrule
\end{tabular}
\end{table}

\begin{table}[!htb]
\centering
\caption{
\nameref{sec:testing:expdesign} rtdet variant solver statistics. 
}
\label{tab:expdesign:rtdettimes}
\small
\sisetup{
table-text-alignment = right,
table-auto-round,
table-figures-integer = 4,
table-figures-decimal = 1,
table-format = 4.1,
add-decimal-zero = false,
add-integer-zero = false,
}
\begin{tabular}{rcrScrScrScrScrS}
\toprule
& \multicolumn{3}{c}{NF} & \multicolumn{6}{c}{EF-exp} & \multicolumn{6}{c}{EF-sec} \\
\cmidrule(lr){2-4} \cmidrule(lr){5-10} \cmidrule(lr){11-16}
& \multicolumn{3}{c}{Hypatia} & \multicolumn{3}{c}{Hypatia} & \multicolumn{3}{c}{MOSEK} & \multicolumn{3}{c}{Hypatia} & \multicolumn{3}{c}{MOSEK} \\
\cmidrule(lr){2-4} \cmidrule(lr){5-7} \cmidrule(lr){8-10} \cmidrule(lr){11-13} \cmidrule(lr){14-16} 
$k$ & {st} & {it} & {time} & {st} & {it} & {time} & {st} & {it} & {time} & {st} & {it} & {time} & {st} & {it} & {time} \\
\midrule
50 & \underline{co} & 25 & 0.3 & \underline{co} & 22 & 4.7 & \underline{co} & 14 & 11. & \underline{co} & 22 & 5.1 & \underline{co} & 11 & 10. \\
100 & \underline{co} & 25 & 0.9 & \underline{co} & 25 & 93. & \underline{co} & 13 & 247. & \underline{co} & 26 & 97. & \underline{co} & 11 & 220. \\
150 & \underline{co} & 26 & 2.6 & \underline{co} & 27 & 696. & \underline{co} & 12 & 1580. & \underline{sp} & 36 & 921. & \underline{co} & 10 & 1432. \\
200 & \underline{co} & 23 & 5.8 & tl & 17 & 1821. & tl & 0 & 1821. & tl & 17 & 1848. & tl & 0 & 1868. \\
300 & \underline{co} & 31 & 31. & sk & $\ast$ & $\ast$ & sk & $\ast$ & $\ast$ & sk & $\ast$ & $\ast$ & sk & $\ast$ & $\ast$ \\
400 & \underline{co} & 29 & 67. & m & $\ast$ & $\ast$ & sk & $\ast$ & $\ast$ & m & $\ast$ & $\ast$ & sk & $\ast$ & $\ast$ \\
500 & \underline{co} & 32 & 152. & sk & $\ast$ & $\ast$ & sk & $\ast$ & $\ast$ & sk & $\ast$ & $\ast$ & sk & $\ast$ & $\ast$ \\
600 & \underline{co} & 33 & 281. & sk & $\ast$ & $\ast$ & sk & $\ast$ & $\ast$ & sk & $\ast$ & $\ast$ & sk & $\ast$ & $\ast$ \\
700 & \underline{co} & 32 & 530. & sk & $\ast$ & $\ast$ & m & $\ast$ & $\ast$ & sk & $\ast$ & $\ast$ & sk & $\ast$ & $\ast$ \\
800 & \underline{co} & 32 & 728. & sk & $\ast$ & $\ast$ & sk & $\ast$ & $\ast$ & sk & $\ast$ & $\ast$ & sk & $\ast$ & $\ast$ \\
900 & \underline{co} & 36 & 1253. & sk & $\ast$ & $\ast$ & sk & $\ast$ & $\ast$ & sk & $\ast$ & $\ast$ & sk & $\ast$ & $\ast$ \\
1000 & \underline{co} & 33 & 1729. & sk & $\ast$ & $\ast$ & sk & $\ast$ & $\ast$ & sk & $\ast$ & $\ast$ & sk & $\ast$ & $\ast$ \\
\bottomrule
\end{tabular}
\end{table}

\subsection{Polynomial minimization}
\label{sec:testing:polymin}

Following \citet{papp2019sum}, we use an interpolant basis weighted SOS dual formulation to find a lower bound for a multivariate polynomial $f$ of maximum degree $2k$ in $m$ variables over the unit hypercube $\mathcal{D} = [-1, 1]^m$.
We let $U = \binom{m+2k}{m}$, $L = \binom{m+k}{m}$, $\tilde{L} = \binom{m+k-1}{m}$.
We select multivariate Chebyshev basis polynomials $g_j, \forall j \in \iin{L}$ of increasing degree up to $k$, and suitable interpolation points $o_u \in \mathcal{D}, \forall u \in \iin{U}$.
To parametrize $\K_{\wsos(P)}^\ast$, we set up the collection of matrices $P$ by evaluating functions of basis polynomials at the points:
\begin{subequations}
\begin{align}
(P_1)_{u, j} & = g_j(o_u) 
& \forall u \in \iin{U}, j \in \iin{L},
\\
(P_{1 + i})_{u, j} & = g_j(o_u) \bigl( 1 - o^2_{u, i} \bigr) 
& \forall i \in \iin{m}, u \in \iin{U}, j \in \iin{\tilde{L}}.
\end{align}
\end{subequations}
Letting $\bar{f} = (f(o_u))_{u \in U}$ be evaluations of $f$ at the points, the conic formulation is:
\begin{subequations}
\begin{align}
\textstyle\min_{\rho \in \bbR^U} \quad \bar{f}' \rho & :
\\
e' \rho & = 1,
\\
\rho & \in \K_{\wsos(P)}^\ast.
\label{eq:polymin:wsos}
\end{align}
\label{eq:polymin}
\end{subequations}
The EF for NF constraint \cref{eq:polymin:wsos} uses $\kpsd$ and is implicit in \cref{eq:wsos:dual}.

To build random instances of \cref{eq:polymin}, we generate $\bar{f}$ (which implicitly defines a polynomial $f$) with independent Gaussian entries, for various values of $m$ and $k$.
Our results are summarized in \cref{tab:polymin}.
Note that $\nu = \bar{\nu}$, $p = \bar{p} = 1$, $n = \bar{n} = q$.
For fixed $m$, the conic constraint dimensions are larger for the EFs and grow much faster for the EFs as the degree $k$ increases.
Hypatia-NF is faster than the EF solvers on all instances with $k > 1$, and solves instances with much higher degrees.
 
\begin{table}[!htb]
\centering
\caption{
\nameref{sec:testing:polymin} formulation and solver statistics. 
}
\label{tab:polymin}
\small
\sisetup{
table-text-alignment = right,
table-auto-round,
table-figures-integer = 4,
table-figures-decimal = 1,
table-format = 4.1,
add-decimal-zero = false,
add-integer-zero = false,
}
\begin{tabular}{rrrrrcrScrScrS}
\toprule
& & \multicolumn{3}{c}{form. stats.} & \multicolumn{3}{c}{Hypatia-NF} & \multicolumn{3}{c}{Hypatia-EF} & \multicolumn{3}{c}{MOSEK-EF} \\
\cmidrule(lr){3-5} \cmidrule(lr){6-8} \cmidrule(lr){9-11} \cmidrule(lr){12-14} 
$m$ & $k$ & $\nu$ & $n$ & $\bar{q}$ & {st} & {it} & {time} & {st} & {it} & {time} & {st} & {it} & {time} \\
\midrule
1 & 100 & 201 & 201 & 10201 & \underline{co} & 12 & 0.1 & \underline{co} & 34 & 1.2 & \underline{co} & 15 & 27. \\
1 & 200 & 401 & 401 & 40401 & \underline{co} & 14 & 0.3 & \underline{co} & 39 & 13. & \underline{co} & 11 & 409. \\
1 & 500 & 1001 & 1001 & 251001 & \underline{co} & 18 & 2.4 & \underline{co} & 57 & 329. & rl & $\ast$ & $\ast$ \\
1 & 1000 & 2001 & 2001 & $\ast$ & \underline{co} & 19 & 11. & m & $\ast$ & $\ast$ & sk & $\ast$ & $\ast$ \\
1 & 2000 & 4001 & 4001 & $\ast$ & \underline{co} & 21 & 73. & sk & $\ast$ & $\ast$ & sk & $\ast$ & $\ast$ \\
1 & 3000 & 6001 & 6001 & $\ast$ & \underline{co} & 24 & 235. & sk & $\ast$ & $\ast$ & sk & $\ast$ & $\ast$ \\
1 & 4000 & 8001 & 8001 & $\ast$ & \underline{co} & 24 & 508. & sk & $\ast$ & $\ast$ & sk & $\ast$ & $\ast$ \\
1 & 5000 & 10001 & 10001 & $\ast$ & \underline{co} & 24 & 916. & sk & $\ast$ & $\ast$ & sk & $\ast$ & $\ast$ \\
2 & 15 & 376 & 496 & 23836 & \underline{co} & 15 & 0.4 & \underline{co} & 21 & 5.0 & \underline{co} & 10 & 87. \\
2 & 30 & 1426 & 1891 & 339946 & \underline{co} & 25 & 10. & \underline{co} & 49 & 751. & rl & $\ast$ & $\ast$ \\
2 & 45 & 3151 & 4186 & $\ast$ & \underline{co} & 22 & 58. & m & $\ast$ & $\ast$ & sk & $\ast$ & $\ast$ \\
2 & 60 & 5551 & 7381 & $\ast$ & \underline{co} & 28 & 300. & sk & $\ast$ & $\ast$ & sk & $\ast$ & $\ast$ \\
2 & 75 & 8626 & 11476 & $\ast$ & \underline{co} & 30 & 1019. & sk & $\ast$ & $\ast$ & sk & $\ast$ & $\ast$ \\
3 & 6 & 252 & 455 & 8358 & \underline{co} & 17 & 0.3 & \underline{co} & 17 & 1.6 & \underline{co} & 9 & 9.1 \\
3 & 9 & 715 & 1330 & 65395 & \underline{co} & 20 & 3.1 & \underline{co} & 24 & 104. & \underline{co} & 9 & 799. \\
3 & 12 & 1547 & 2925 & 303030 & \underline{co} & 23 & 20. & \underline{co} & 33 & 1775. & rl & $\ast$ & $\ast$ \\
3 & 15 & 2856 & 5456 & $\ast$ & \underline{co} & 26 & 89. & m & $\ast$ & $\ast$ & sk & $\ast$ & $\ast$ \\
3 & 18 & 4750 & 9139 & $\ast$ & er & 34 & 1340. & sk & $\ast$ & $\ast$ & sk & $\ast$ & $\ast$ \\
4 & 4 & 210 & 495 & 5005 & \underline{co} & 18 & 0.4 & \underline{co} & 16 & 1.7 & \underline{co} & 8 & 3.9 \\
4 & 6 & 714 & 1820 & 54159 & \underline{co} & 15 & 4.8 & \underline{co} & 18 & 222. & \underline{co} & 10 & 579. \\
4 & 8 & 1815 & 4845 & $\ast$ & \underline{co} & 20 & 63. & m & $\ast$ & $\ast$ & m & $\ast$ & $\ast$ \\
4 & 10 & 3861 & 10626 & $\ast$ & \underline{co} & 22 & 458. & sk & $\ast$ & $\ast$ & sk & $\ast$ & $\ast$ \\
8 & 2 & 117 & 495 & 1395 & \underline{co} & 26 & 0.5 & \underline{co} & 21 & 0.7 & \underline{co} & 11 & 0.9 \\
8 & 3 & 525 & 3003 & 21975 & \underline{co} & 18 & 15. & \underline{co} & 16 & 148. & \underline{co} & 8 & 125. \\
8 & 4 & 1815 & 12870 & $\ast$ & \underline{co} & 27 & 633. & m & $\ast$ & $\ast$ & m & $\ast$ & $\ast$ \\
16 & 1 & 33 & 153 & 169 & \underline{co} & 13 & 0.1 & \underline{co} & 12 & 0.7 & \underline{co} & 7 & 0.0 \\
16 & 2 & 425 & 4845 & 14229 & \underline{co} & 27 & 86. & \underline{co} & 22 & 174. & \underline{sp} & 10 & 192. \\
32 & 1 & 65 & 561 & 593 & \underline{co} & 15 & 0.7 & \underline{co} & 12 & 1.0 & \underline{co} & 7 & 0.2 \\
64 & 1 & 129 & 2145 & 2209 & \underline{co} & 15 & 14. & \underline{co} & 12 & 3.1 & \underline{co} & 9 & 3.0 \\
\bottomrule
\end{tabular}
\end{table}
\subsection{Smooth density estimation}
\label{sec:testing:densityest}

$\bbR_{m,2k}[x]$ is the ring of polynomials of maximum degree $2k$ in $m$ variables \citep{papp2019sum}.
We seek a polynomial density function $f \in \bbR_{m,2k}[x]$ over the domain $\mathcal{D} = [-1, 1]^m$ that maximizes the log-likelihood of $N$ given observations $z_i \in \mathcal{D}, \forall i \in \iin{N}$ (compare to \citet[Section 4.3]{papp2014shape}).
For $f$ to be a valid density it must be nonnegative on $\mathcal{D}$ and integrate to one over $\mathcal{D}$, so we aim to solve:
\begin{subequations}
\begin{align}
\textstyle\max_{f \in \bbR_{m,2k}[x]} \quad \tsum{i \in \iin{N}} \log(f(z_i)) & : 
\\
\textstyle\int_{\mathcal{D}} f(x) \, d x & = 1,
\label{eq:densityest:infinite:integral}
\\
f(x) & \geq 0 
\quad \forall x \in \mathcal{D}.
\end{align}
\label{eq:densityest:infinite}
\end{subequations}
To find a feasible solution for \cref{eq:densityest:infinite}, we build an SOS formulation.
We obtain interpolation points and matrices $P$ parametrizing $\K_{\wsos(P)}$, using the techniques from \cref{sec:testing:polymin}.
From the interpolation points and the domain $\mathcal{D}$, we compute a vector of quadrature weights $\mu \in \bbR^U$.
We compute a matrix $B \in \bbR^{N \times U}$ by evaluating the $U$ Lagrange basis polynomials corresponding to the interpolation points (see \citet{papp2019sum}) at the $N$ observations.
Letting variable $\rho$ represent the coefficients on the Lagrange basis, the conic formulation is:
\begin{subequations}
\begin{align}
\textstyle\max_{\psi \in \bbR, \rho \in \bbR^U} \quad \psi & : 
\\
\mu' \rho & = 1,
\\
(\psi, 1, B \rho) & \in \klog,
\label{eq:densityest:log}
\\
\rho & \in \K_{\wsos(P)}.
\label{eq:densityest:wsos}
\end{align}
\label{eq:densityest}
\end{subequations}
The EFs for NF constraints \cref{eq:densityest:log,eq:densityest:wsos} follow \cref{eq:cones:log,eq:wsos:prim}.

To build random instances of \cref{eq:densityest} for various values of $m$ and $k$, we generate $N = 500$ independent uniform samples in $[-1, 1]^m$ for $z_i \in \mathcal{D}, \forall i \in \iin{N}$.
As our method for computing $\mu$ is numerically unstable for larger $m$, we only use $m \leq 16$.
Our results are summarized in \cref{tab:densityest}.
Note that $\bar{\nu} = 999 + \nu$, $p = 1$, $\bar{p} = n$, $q = 501 + n$, $\bar{q} = 1001 + \bar{n} - n$.
All dimensions are larger for the EFs than for the NFs.
Hypatia-NF is faster than the EF solvers and solves instances with much higher degrees. 

\begin{table}[!htb]
\centering
\caption{
\nameref{sec:testing:densityest}. 
}
\label{tab:densityest}
\small
\sisetup{
table-text-alignment = right,
table-auto-round,
table-figures-integer = 4,
table-figures-decimal = 1,
table-format = 4.1,
add-decimal-zero = false,
add-integer-zero = false,
}
\begin{tabular}{rrrrrcrScrScrS}
\toprule
& & \multicolumn{3}{c}{dimensions} & \multicolumn{3}{c}{Hypatia-NF} & \multicolumn{3}{c}{Hypatia-EF} & \multicolumn{3}{c}{MOSEK-EF} \\
\cmidrule(lr){3-5} \cmidrule(lr){6-8} \cmidrule(lr){9-11} \cmidrule(lr){12-14} 
$m$ & $2k$ & $\nu$ & $n$ & $\bar{n}$ & {st} & {it} & {time} & {st} & {it} & {time} & {st} & {it} & {time} \\
\midrule
1 & 250 & 753 & 252 & 16628 & \underline{co} & 35 & 0.2 & sp & 43 & 1040. & \underline{co} & 25 & 112. \\
1 & 500 & 1003 & 502 & 64003 & \underline{co} & 42 & 1.1 & rl & $\ast$ & $\ast$ & tl & 18 & 1814. \\
1 & 1000 & 1503 & 1002 & $\ast$ & \underline{co} & 41 & 5.1 & m & $\ast$ & $\ast$ & m & $\ast$ & $\ast$ \\
1 & 2000 & 2503 & 2002 & $\ast$ & \underline{co} & 56 & 23. & sk & $\ast$ & $\ast$ & sk & $\ast$ & $\ast$ \\
1 & 4000 & 4503 & 4002 & $\ast$ & \underline{co} & 82 & 185. & sk & $\ast$ & $\ast$ & sk & $\ast$ & $\ast$ \\
1 & 6000 & 6503 & 6002 & $\ast$ & \underline{co} & 106 & 663. & sk & $\ast$ & $\ast$ & sk & $\ast$ & $\ast$ \\
2 & 20 & 678 & 232 & 6023 & \underline{co} & 50 & 0.3 & \underline{co} & 29 & 73. & \underline{co} & 19 & 7.4 \\
2 & 40 & 1153 & 862 & 72468 & \underline{co} & 34 & 2.9 & rl & $\ast$ & $\ast$ & sp & 22 & 1522. \\
2 & 60 & 1928 & 1892 & $\ast$ & \underline{co} & 36 & 11. & m & $\ast$ & $\ast$ & m & $\ast$ & $\ast$ \\
2 & 80 & 3003 & 3322 & $\ast$ & \underline{co} & 53 & 64. & sk & $\ast$ & $\ast$ & sk & $\ast$ & $\ast$ \\
2 & 100 & 4378 & 5152 & $\ast$ & \underline{co} & 64 & 247. & sk & $\ast$ & $\ast$ & sk & $\ast$ & $\ast$ \\
3 & 12 & 754 & 456 & 9314 & \underline{co} & 57 & 0.9 & \underline{co} & 25 & 293. & sp & 19 & 20. \\
3 & 18 & 1217 & 1331 & 67226 & \underline{co} & 55 & 6.6 & rl & $\ast$ & $\ast$ & sp & 17 & 1216. \\
3 & 24 & 2049 & 2926 & $\ast$ & \underline{co} & 46 & 35. & m & $\ast$ & $\ast$ & m & $\ast$ & $\ast$ \\
3 & 30 & 3358 & 5457 & $\ast$ & \underline{co} & 63 & 348. & sk & $\ast$ & $\ast$ & sk & $\ast$ & $\ast$ \\
4 & 8 & 712 & 496 & 6001 & \underline{co} & 57 & 1.0 & \underline{co} & 25 & 133. & sp & 22 & 12. \\
4 & 12 & 1216 & 1821 & 56480 & \underline{co} & 72 & 16. & tl & $\ast$ & $\ast$ & sp & 20 & 934. \\
4 & 16 & 2317 & 4846 & $\ast$ & \underline{co} & 70 & 192. & m & $\ast$ & $\ast$ & m & $\ast$ & $\ast$ \\
8 & 4 & 619 & 496 & 2391 & \underline{co} & 96 & 2.1 & \underline{co} & 30 & 9.8 & sp & 17 & 2.1 \\
8 & 6 & 1027 & 3004 & 25479 & \underline{co} & 90 & 62. & tl & $\ast$ & $\ast$ & sp & 17 & 379. \\
\bottomrule
\end{tabular}
\end{table}

\subsection{Shape constrained regression}
\label{sec:testing:shapeconregr}

A common type of shape constraint imposes monotonicity or convexity of a polynomial over a basic semialgebraic set \citep[Section 6]{hall2019engineering}.
Given an $m$-dimensional feature variable $z$ and a scalar response variable $g$, we aim to fit a polynomial $f \in \bbR_{m,2k}[x]$ that is convex over $\mathcal{D} = [-1, 1]^m$ to $N$ given observations $(z_i, g_i)_{i \in \iin{N}}$ with $z_i \in \mathcal{D}, \forall i \in \iin{N}$:
\begin{subequations}
\begin{align}
\textstyle\min_{f \in \bbR_{m,2k}[x]} \quad 
\tsum{i \in \iin{N}} (g_i - f(z_i))^2 & :
\\
y' ( \nabla^2 f(x) ) y & \geq 0 
\quad \forall x \in \mathcal{D}, y \in \bbR^m.
\label{eq:shapeconregr:infinite:matwsos}
\end{align}
\label{eq:shapeconregr:infinite}
\end{subequations}
Constraint \cref{eq:shapeconregr:infinite:matwsos} ensures the Hessian matrix $\nabla^2 f(x)$ of polynomials is PSD at every point $x \in \mathcal{D}$, which is equivalent to convexity of $f$ over $\mathcal{D}$.
To find a feasible solution for \cref{eq:shapeconregr:infinite}, we build an SOS formulation.
The polynomial variable, represented in an interpolant basis with the optimization variable $\rho \in \bbR^U$, has degree $2k$ and $U = \binom{m + 2k}{m}$ coefficients. 
Each polynomial entry of $\nabla^2 f(x)$ has degree $2k-2$ and $\bar{U} = \binom{m + 2k-2}{m}$ coefficients.
Following the descriptions in \crefrange{sec:testing:polymin}{sec:testing:densityest}, we obtain interpolation points and a Lagrange polynomial basis for these $U$-dimensional and $\bar{U}$-dimensional spaces, and we define the matrix $B \in \bbR^{N \times U}$ containing evaluations of the $U$-dimensional Lagrange basis at the $N$ feature observations.
Finally, we let $F \in \bbR^{\sdim(m) \bar{U} \times U}$ be such that $F \rho$ is a vectorization of the tensor $H \in \bbR^{m \times m \times \bar{U}}$ (scaled to account for symmetry) with $H_{a,b,u}$ equal to the $u$th coefficient of the $(a,b)$th polynomial in $\nabla^2 f(x)$ for $a,b \in \iin{m}$ and $u \in \iin{\bar{U}}$.
This yields the formulation:
\begin{subequations}
\begin{align}
\textstyle\min_{\psi \in \bbR, \rho \in \bbR^U} \quad \psi & :
\\
(\psi, g - B \rho) & \in \kltwo,
\label{eq:shapeconregr:kl2}
\\
F \rho & \in \K_{\matwsos(P)}.
\label{eq:shapeconregr:matwsos}
\end{align}
\label{eq:shapeconregr}
\end{subequations}
Note that for $N > U$, we use a QR factorization to reduce the dimension of $\kltwo$ in \cref{eq:shapeconregr:kl2} from $1 + N$ to $2 + U$.\footnote{
Let $[-B \ g] = Q R$, where $Q \in \bbR^{N \times (U + 1)}$ has orthonormal columns and $R \in \bbR^{(U + 1) \times (U + 1)}$ is upper triangular. 
Then $(\psi, g - B \rho) \in \kltwo$ if and only if $(\psi, R (\rho, 1)) \in \kltwo$. 
}
The EF for NF constraint \cref{eq:shapeconregr:matwsos} follows \cref{eq:matwsos:prim}.

To build random instances of \cref{eq:shapeconregr} for various values of $m$ and $k$, we generate $N = \lceil 1.1 U \rceil$ independent observations with $z_i$ sampled uniformly from $\mathcal{D}$ and $g_i = \exp ( \lVert z \rVert^2 / m) - 1 + \varepsilon_i$, where $\varepsilon_i$ is a Gaussian sample yielding a signal to noise ratio of $10$, for all $i \in \iin{N}$.
We exclude the case $m = 1$, since $\K_{\wsos(P)}$ could be used in place of $\K_{\matwsos(P)}$.
Our results are summarized in \cref{tab:shapeconregr}.
Note that $\nu = \bar{\nu}$, $p = 0$, $\bar{p} = q - n - 1$, $\bar{q} = \bar{n} + 1$.
All dimensions are larger for the EFs.
The instances are numerically challenging, and MOSEK-EF often encounters slow progress.
Hypatia-NF is faster than the EF solvers and solves instances with much higher degrees.

\begin{table}[!htb]
\centering
\caption{
\nameref{sec:testing:shapeconregr} formulation and solver statistics. 
}
\label{tab:shapeconregr}
\small
\sisetup{
table-text-alignment = right,
table-auto-round,
table-figures-integer = 4,
table-figures-decimal = 1,
table-format = 4.1,
add-decimal-zero = false,
add-integer-zero = false,
}
\begin{tabular}{rrrrrrcrScrScrS}
\toprule
& & \multicolumn{4}{c}{form. stats.} & \multicolumn{3}{c}{Hypatia-NF} & \multicolumn{3}{c}{Hypatia-EF} & \multicolumn{3}{c}{MOSEK-EF} \\
\cmidrule(lr){3-6} \cmidrule(lr){7-9} \cmidrule(lr){10-12} \cmidrule(lr){13-15} 
$m$ & $2k$ & $\nu$ & $n$ & $\bar{n}$ & $q$ & {st} & {it} & {time} & {st} & {it} & {time} & {st} & {it} & {time} \\
\midrule
2 & 10 & 72 & 67 & 952 & 203 & \underline{co} & 18 & 0.0 & \underline{co} & 24 & 0.6 & \underline{sp} & 19 & 0.4 \\
2 & 20 & 292 & 232 & 14527 & 803 & \underline{co} & 37 & 0.5 & \underline{sp} & 68 & 1348. & sp & 20 & 60. \\
2 & 30 & 662 & 497 & 73727 & 1803 & \underline{co} & 58 & 5.6 & rl & $\ast$ & $\ast$ & tl & 12 & 1845. \\
2 & 40 & 1182 & 862 & $\ast$ & 3203 & \underline{co} & 84 & 34. & m & $\ast$ & $\ast$ & sk & $\ast$ & $\ast$ \\
2 & 50 & 1852 & 1327 & $\ast$ & 5003 & er & 52 & 161. & sk & $\ast$ & $\ast$ & sk & $\ast$ & $\ast$ \\
2 & 60 & 2672 & 1892 & $\ast$ & 7203 & \underline{er} & 120 & 939. & sk & $\ast$ & $\ast$ & sk & $\ast$ & $\ast$ \\
3 & 8 & 152 & 166 & 3391 & 671 & \underline{co} & 19 & 0.1 & \underline{co} & 27 & 14. & \underline{sp} & 23 & 4.3 \\
3 & 12 & 485 & 456 & 31347 & 2173 & \underline{co} & 38 & 3.5 & tl & $\ast$ & $\ast$ & sp & 15 & 244. \\
3 & 16 & 1118 & 970 & 161584 & 5051 & \underline{co} & 61 & 44. & m & $\ast$ & $\ast$ & rl & $\ast$ & $\ast$ \\
3 & 20 & 2147 & 1772 & $\ast$ & 9753 & \underline{co} & 87 & 325. & sk & $\ast$ & $\ast$ & sk & $\ast$ & $\ast$ \\
3 & 24 & 3668 & 2926 & $\ast$ & 16727 & \underline{co} & 111 & 1605. & sk & $\ast$ & $\ast$ & sk & $\ast$ & $\ast$ \\
4 & 6 & 142 & 211 & 2881 & 912 & \underline{co} & 17 & 0.4 & \underline{co} & 23 & 7.8 & \underline{sp} & 18 & 3.7 \\
4 & 8 & 382 & 496 & 17686 & 2597 & \underline{co} & 24 & 3.1 & \underline{co} & 38 & 1621. & sp & 14 & 94. \\
4 & 10 & 842 & 1002 & 79822 & 5953 & \underline{co} & 38 & 35. & rl & $\ast$ & $\ast$ & tl & 12 & 2068. \\
4 & 12 & 1626 & 1821 & $\ast$ & 11832 & \underline{co} & 58 & 283. & m & $\ast$ & $\ast$ & sk & $\ast$ & $\ast$ \\
4 & 14 & 2858 & 3061 & $\ast$ & 21262 & \underline{co} & 72 & 1430. & sk & $\ast$ & $\ast$ & sk & $\ast$ & $\ast$ \\
6 & 4 & 80 & 211 & 1240 & 800 & \underline{co} & 13 & 0.3 & \underline{co} & 15 & 0.7 & \underline{co} & 9 & 0.6 \\
6 & 6 & 422 & 925 & 20539 & 5336 & \underline{co} & 22 & 15. & \underline{co} & 32 & 1630. & \underline{sp} & 17 & 260. \\
6 & 8 & 1514 & 3004 & 215440 & 22409 & \underline{co} & 29 & 638. & m & $\ast$ & $\ast$ & rl & $\ast$ & $\ast$ \\
8 & 4 & 138 & 496 & 3412 & 2117 & \underline{co} & 19 & 2.1 & \underline{co} & 20 & 8.0 & \underline{co} & 11 & 4.1 \\
8 & 6 & 938 & 3004 & 89008 & 20825 & \underline{co} & 31 & 515. & rl & $\ast$ & $\ast$ & rl & $\ast$ & $\ast$ \\
10 & 4 & 212 & 1002 & 7657 & 4633 & \underline{co} & 26 & 13. & \underline{co} & 24 & 87. & \underline{co} & 13 & 25. \\
12 & 4 & 302 & 1821 & 15003 & 8920 & \underline{co} & 29 & 73. & \underline{co} & 28 & 504. & \underline{co} & 11 & 125. \\
14 & 4 & 408 & 3061 & 26686 & 15662 & \underline{co} & 33 & 346. & tl & 0 & 1884. & \underline{sp} & 21 & 767. \\
\bottomrule
\end{tabular}
\end{table}

\section{Conclusions}
\label{sec:conclusion}

Although many convex problems are representable with conic EFs using the small number of standard cones currently recognized by some advanced conic solvers, these formulations can be much larger and more complex than NFs with exotic cones.
In \cref{sec:cones}, we describe some of Hypatia's predefined exotic cones and analyze general techniques for constructing EFs from NFs that use these cones.
For several example problems, we propose NFs and generate instances of a wide range of sizes.
Across these instances, we observe much higher empirical dimensions (variable, equality, and conic constraint dimensions in the conic general form \cref{eq:prim}) for the EFs than for the NFs.
We demonstrate significant computational advantages from solving the NFs with Hypatia compared to solving the EFs with either Hypatia or MOSEK 9, especially in terms of solve time and memory usage.
We also observe that the NFs are typically faster and less memory-intensive to generate using JuMP.

Our results suggest that when there exists an NF that is significantly smaller than any EF, it is probably worth trying to solve the NF with Hypatia.
In deciding whether to formulate an NF or an EF, it can be helpful to examine our summary in \cref{tab:natvsext} of computational properties for NFs and EFs of exotic cone constraints.
For spectral and nuclear norm constraints, when the matrix ($W \in \bbR^{d_1 \times d_2}$) has many more columns than rows ($d_2 \gg d_1$), the dimensions look relatively more favorable for the NF.
For SOS and SOS matrix constraints, the dimensions grow much more slowly for the NF as the polynomial degree increases. 
Sometimes the modeler has to choose between different EFs.
For our matrix completion problem and experiment design root-determinant variant, we compare two EFs for the geometric mean cone and find that Hypatia performs better with the exponential cone EF (EF-exp) and MOSEK performs better with the second order cone EF (EF-sec).

If the modeler has an NF that uses a proper cone not already defined in Hypatia, the user can add support for the cone through Hypatia's generic cone interface.
It may require some effort to make the cone oracles as efficient and numerically stable as possible.
However, we have already predefined over two dozen exotic cone types with tractable oracles in Hypatia \citep{coey2021performance,coey2021hypatia}, many of which have multiple variants (such as real or complex flavors).
We use these cones to model hundreds of formulations in our over three dozen examples available in Hypatia's examples folder. 
Our experience with these examples suggests that NFs tend to be more convenient for modeling and interpreting conic certificates.

\bibliographystyle{abbrvnat}
\bibliography{refs}

\begin{thebibliography}{43}
\providecommand{\natexlab}[1]{#1}
\providecommand{\url}[1]{\texttt{#1}}
\expandafter\ifx\csname urlstyle\endcsname\relax
  \providecommand{\doi}[1]{doi: #1}\else
  \providecommand{\doi}{doi: \begingroup \urlstyle{rm}\Url}\fi

\bibitem[Agrawal et~al.(2019)Agrawal, Diamond, and
  Boyd]{agrawal2019disciplined}
A.~Agrawal, S.~Diamond, and S.~Boyd.
\newblock Disciplined geometric programming.
\newblock \emph{Optimization Letters}, 13\penalty0 (5):\penalty0 961--976,
  2019.

\bibitem[Andersen et~al.(2011)Andersen, Dahl, Liu, Vandenberghe, Sra, Nowozin,
  and Wright]{andersen2011interior}
M.~Andersen, J.~Dahl, Z.~Liu, L.~Vandenberghe, S.~Sra, S.~Nowozin, and
  S.~Wright.
\newblock Interior-point methods for large-scale cone programming.
\newblock \emph{Optimization for Machine Learning}, 5583, 2011.

\bibitem[Ben-Tal and Nemirovski(2001)]{ben2001lectures}
A.~Ben-Tal and A.~Nemirovski.
\newblock \emph{Lectures on modern convex optimization: analysis, algorithms,
  and engineering applications}, volume~2.
\newblock {SIAM}, 2001.

\bibitem[Bezanson et~al.(2017)Bezanson, Edelman, Karpinski, and
  Shah]{bezanson2017julia}
J.~Bezanson, A.~Edelman, S.~Karpinski, and V.~B. Shah.
\newblock Julia: A fresh approach to numerical computing.
\newblock \emph{SIAM Review}, 59\penalty0 (1):\penalty0 65--98, 2017.

\bibitem[Borchers(1999)]{borchers1999csdp}
B.~Borchers.
\newblock {CSDP}, a {C} library for semidefinite programming.
\newblock \emph{Optimization Methods and Software}, 11\penalty0 (1-4):\penalty0
  613--623, 1999.

\bibitem[Boyd and Vandenberghe(2004)]{boyd2004convex}
S.~Boyd and L.~Vandenberghe.
\newblock \emph{Convex optimization}.
\newblock Cambridge University Press, 2004.

\bibitem[Coey et~al.(2020)Coey, Lubin, and Vielma]{coey2020outer}
C.~Coey, M.~Lubin, and J.~P. Vielma.
\newblock Outer approximation with conic certificates for mixed-integer convex
  problems.
\newblock \emph{Mathematical Programming Computation}, 12:\penalty0 249--293,
  2020.

\bibitem[Coey et~al.(2021{\natexlab{a}})Coey, Kapelevich, and
  Vielma]{coey2021hypatia}
C.~Coey, L.~Kapelevich, and J.~P. Vielma.
\newblock Hypatia cones reference, 2021{\natexlab{a}}.
\newblock URL \url{https://github.com/chriscoey/Hypatia.jl/wiki/}.
\newblock Online; accessed 1-June-2021.

\bibitem[Coey et~al.(2021{\natexlab{b}})Coey, Kapelevich, and
  Vielma]{coey2021hypatiab}
C.~Coey, L.~Kapelevich, and J.~P. Vielma.
\newblock \emph{Hypatia documentation}, 2021{\natexlab{b}}.
\newblock URL \url{https://chriscoey.github.io/Hypatia.jl/dev/}.
\newblock Online; accessed 7-June-2021.

\bibitem[Coey et~al.(2021{\natexlab{c}})Coey, Kapelevich, and
  Vielma]{coey2021performance}
C.~Coey, L.~Kapelevich, and J.~P. Vielma.
\newblock Performance enhancements for a generic conic interior point
  algorithm.
\newblock \emph{arXiv preprint arXiv:?????}, 2021{\natexlab{c}}.

\bibitem[Coey et~al.(2021{\natexlab{d}})Coey, Kapelevich, and
  Vielma]{coey2021self}
C.~Coey, L.~Kapelevich, and J.~P. Vielma.
\newblock Conic optimization over epigraphs of spectral functions on cones of
  squares.
\newblock \emph{arXiv preprint arXiv:2103.04104}, 2021{\natexlab{d}}.

\bibitem[Dahl and Andersen(2021)]{dahl2021primal}
J.~Dahl and E.~D. Andersen.
\newblock A primal-dual interior-point algorithm for nonsymmetric
  exponential-cone optimization.
\newblock \emph{Mathematical Programming}, pages 1--30, 2021.

\bibitem[Diamond and Boyd(2016)]{diamond2016cvxpy}
S.~Diamond and S.~Boyd.
\newblock {CVXPY:} a {Python}-embedded modeling language for convex
  optimization.
\newblock \emph{The Journal of Machine Learning Research}, 17\penalty0
  (1):\penalty0 2909--2913, 2016.

\bibitem[Domahidi et~al.(2013)Domahidi, Chu, and Boyd]{domahidi2013ecos}
A.~Domahidi, E.~Chu, and S.~Boyd.
\newblock {ECOS:} an {SOCP} solver for embedded systems.
\newblock In \emph{2013 European Control Conference (ECC)}, pages 3071--3076.
  IEEE, 2013.

\bibitem[Dunning et~al.(2017)Dunning, Huchette, and Lubin]{dunning2017jump}
I.~Dunning, J.~Huchette, and M.~Lubin.
\newblock {JuMP:} a modeling language for mathematical optimization.
\newblock \emph{SIAM Review}, 59\penalty0 (2):\penalty0 295--320, 2017.

\bibitem[Grant and Boyd(2014)]{grant2014cvx}
M.~Grant and S.~Boyd.
\newblock {CVX}: {MATLAB} software for disciplined convex programming, version
  2.1, 2014.

\bibitem[G{\"u}ler(1996)]{guler1996barrier}
O.~G{\"u}ler.
\newblock Barrier functions in interior point methods.
\newblock \emph{Mathematics of Operations Research}, 21\penalty0 (4):\penalty0
  860--885, 1996.

\bibitem[Hall(2019)]{hall2019engineering}
G.~Hall.
\newblock Engineering and business applications of sum of squares polynomials.
\newblock \emph{arXiv preprint arXiv:1906.07961}, 2019.

\bibitem[Kapelevich et~al.(2021)Kapelevich, Coey, and
  Vielma]{kapelevich2021sum}
L.~Kapelevich, C.~Coey, and J.~P. Vielma.
\newblock Sum of squares generalizations for conic sets.
\newblock \emph{arXiv preprint arXiv:2103.11499}, 2021.

\bibitem[Karimi and Tun{\c{c}}el(2019)]{karimi2019domain}
M.~Karimi and L.~Tun{\c{c}}el.
\newblock Domain-driven solver ({DDS}): a {MATLAB}-based software package for
  convex optimization problems in domain-driven form.
\newblock \emph{arXiv preprint arXiv:1908.03075}, 2019.

\bibitem[Karimi and Tun{\c{c}}el(2020)]{karimi2020primal}
M.~Karimi and L.~Tun{\c{c}}el.
\newblock Primal--dual interior-point methods for domain-driven formulations.
\newblock \emph{Mathematics of Operations Research}, 45\penalty0 (2):\penalty0
  591--621, 2020.

\bibitem[Legat et~al.(2020)Legat, Dowson, Garcia, and
  Lubin]{legat2020mathoptinterface}
B.~Legat, O.~Dowson, J.~D. Garcia, and M.~Lubin.
\newblock {MathOptInterface}: a data structure for mathematical optimization
  problems.
\newblock \emph{arXiv preprint arXiv:2002.03447}, 2020.

\bibitem[Lubin et~al.(2016)Lubin, Yamangil, Bent, and
  Vielma]{lubin2016extended}
M.~Lubin, E.~Yamangil, R.~Bent, and J.~P. Vielma.
\newblock Extended formulations in mixed-integer convex programming.
\newblock In Q.~Louveaux and M.~Skutella, editors, \emph{Integer Programming
  and Combinatorial Optimization: 18\textsuperscript{th} International
  Conference, IPCO 2016, Li\`{e}ge, Belgium, June 1-3, 2016, Proceedings},
  pages 102--113. Springer International Publishing, 2016.
\newblock ISBN 978-3-319-33461-5.
\newblock \doi{10.1007/978-3-319-33461-5_9}.

\bibitem[{MOSEK ApS}(2020{\natexlab{a}})]{aps2020mosek}
{MOSEK ApS}.
\newblock {MOSEK} fusion {API} for {Python}, 2020{\natexlab{a}}.
\newblock URL \url{https://docs.mosek.com/9.1/pythonfusion/index.html}.

\bibitem[{MOSEK ApS}(2020{\natexlab{b}})]{mosek2020modeling}
{MOSEK ApS}.
\newblock {Modeling Cookbook} revision 3.2.1, 2020{\natexlab{b}}.
\newblock URL \url{https://docs.mosek.com/modeling-cookbook/index.html}.

\bibitem[Nesterov(2006)]{nesterov2006constructing}
Y.~Nesterov.
\newblock Constructing self-concordant barriers for convex cones.
\newblock \emph{CORE discussion paper}, 2006.

\bibitem[Nesterov(2012)]{nesterov2012towards}
Y.~Nesterov.
\newblock Towards non-symmetric conic optimization.
\newblock \emph{Optimization Methods and Software}, 27\penalty0 (4-5):\penalty0
  893--917, 2012.

\bibitem[Nesterov and Nemirovskii(1994)]{nesterov1994interior}
Y.~Nesterov and A.~Nemirovskii.
\newblock \emph{Interior-point polynomial algorithms in convex programming}.
\newblock Studies in Applied Mathematics. Society for Industrial and Applied
  Mathematics, 1994.

\bibitem[Nesterov and Todd(1997)]{nesterov1997self}
Y.~E. Nesterov and M.~J. Todd.
\newblock Self-scaled barriers and interior-point methods for convex
  programming.
\newblock \emph{Mathematics of Operations Research}, 22\penalty0 (1):\penalty0
  1--42, 1997.

\bibitem[Nesterov et~al.(1996)Nesterov, Todd, and Ye]{nesterov1996infeasible}
Y.~E. Nesterov, M.~J. Todd, and Y.~Ye.
\newblock Infeasible-start primal-dual methods and infeasibility detectors for
  nonlinear programming problems.
\newblock Technical report, Cornell University Operations Research and
  Industrial Engineering, 1996.

\bibitem[O’Donoghue et~al.(2016)O’Donoghue, Chu, Parikh, and
  Boyd]{o2016conic}
B.~O’Donoghue, E.~Chu, N.~Parikh, and S.~Boyd.
\newblock Conic optimization via operator splitting and homogeneous self-dual
  embedding.
\newblock \emph{Journal of Optimization Theory and Applications}, 169\penalty0
  (3):\penalty0 1042--1068, 2016.

\bibitem[Papp and Alizadeh(2014)]{papp2014shape}
D.~Papp and F.~Alizadeh.
\newblock Shape-constrained estimation using nonnegative splines.
\newblock \emph{Journal of Computational and Graphical Statistics}, 23\penalty0
  (1):\penalty0 211--231, 2014.

\bibitem[Papp and Y{\i}ld{\i}z(2017)]{papp2017homogeneous}
D.~Papp and S.~Y{\i}ld{\i}z.
\newblock On ``{A} homogeneous interior-point algorithm for non-symmetric
  convex conic optimization''.
\newblock \emph{arXiv preprint arXiv:1712.00492}, 2017.

\bibitem[Papp and Yildiz(2019)]{papp2019sum}
D.~Papp and S.~Yildiz.
\newblock Sum-of-squares optimization without semidefinite programming.
\newblock \emph{SIAM Journal on Optimization}, 29\penalty0 (1):\penalty0
  822--851, 2019.

\bibitem[Papp and Y{\i}ld{\i}z(2021)]{papp2021alfonso}
D.~Papp and S.~Y{\i}ld{\i}z.
\newblock alfonso: Matlab package for nonsymmetric conic optimization.
\newblock \emph{arXiv preprint arXiv:2101.04274}, 2021.

\bibitem[Permenter et~al.(2017)Permenter, Friberg, and
  Andersen]{permenter2017solving}
F.~Permenter, H.~A. Friberg, and E.~D. Andersen.
\newblock Solving conic optimization problems via self-dual embedding and
  facial reduction: a unified approach.
\newblock \emph{SIAM Journal on Optimization}, 27\penalty0 (3):\penalty0
  1257--1282, 2017.

\bibitem[Recht et~al.(2010)Recht, Fazel, and Parrilo]{recht2010guaranteed}
B.~Recht, M.~Fazel, and P.~A. Parrilo.
\newblock Guaranteed minimum-rank solutions of linear matrix equations via
  nuclear norm minimization.
\newblock \emph{SIAM Review}, 52\penalty0 (3):\penalty0 471--501, 2010.

\bibitem[Serrano(2015)]{serrano2015algorithms}
S.~A. Serrano.
\newblock \emph{Algorithms for unsymmetric cone optimization and an
  implementation for problems with the exponential cone}.
\newblock PhD thesis, Stanford University, 2015.

\bibitem[Skajaa and Ye(2015)]{skajaa2015homogeneous}
A.~Skajaa and Y.~Ye.
\newblock A homogeneous interior-point algorithm for nonsymmetric convex conic
  optimization.
\newblock \emph{Mathematical Programming}, 150\penalty0 (2):\penalty0 391--422,
  2015.

\bibitem[Udell et~al.(2014)Udell, Mohan, Zeng, Hong, Diamond, and
  Boyd]{udell2014convex}
M.~Udell, K.~Mohan, D.~Zeng, J.~Hong, S.~Diamond, and S.~Boyd.
\newblock Convex optimization in {Julia}.
\newblock In \emph{Proceedings of the 1st First Workshop for High Performance
  Technical Computing in Dynamic Languages}, pages 18--28. IEEE Press, 2014.

\bibitem[Vandenberghe(2010)]{vandenberghe2010cvxopt}
L.~Vandenberghe.
\newblock The {CVXOPT} linear and quadratic cone program solvers, 2010.
\newblock URL
  \url{https://www.seas.ucla.edu/~vandenbe/publications/coneprog.pdf}.

\bibitem[Yamashita et~al.(2003)Yamashita, Fujisawa, and
  Kojima]{yamashita2003implementation}
M.~Yamashita, K.~Fujisawa, and M.~Kojima.
\newblock Implementation and evaluation of {SDPA} 6.0 (semidefinite programming
  algorithm 6.0).
\newblock \emph{Optimization Methods and Software}, 18\penalty0 (4):\penalty0
  491--505, 2003.

\bibitem[Yang et~al.(2016)Yang, Luo, Qian, Tai, Zhang, and Xu]{yang2016nuclear}
J.~Yang, L.~Luo, J.~Qian, Y.~Tai, F.~Zhang, and Y.~Xu.
\newblock Nuclear norm based matrix regression with applications to face
  recognition with occlusion and illumination changes.
\newblock \emph{IEEE Transactions on Pattern Analysis and Machine
  Intelligence}, 39\penalty0 (1):\penalty0 156--171, 2016.

\end{thebibliography}

\end{document}